# FLUCTUATION THEORY OF CONNECTIVITIES FOR SUBCRITICAL RANDOM CLUSTER MODELS

By Massimo Campanino, Dmitry Ioffe [1] and Yvan Velenik [2]

*Università di Bologna, Technion and Université de Genève*

We develop a fluctuation theory of connectivities for subcritical random cluster models. The theory is based on a comprehensive nonperturbative probabilistic description of long connected clusters in terms of essentially one-dimensional chains of irreducible objects. Statistics of local observables, for example, displacement, over such chains obey classical limit laws, and our construction leads to an effective random walk representation of percolation clusters.

The results include a derivation of a sharp Ornstein–Zernike type asymptotic formula for two point functions, a proof of analyticity and strict convexity of inverse correlation length and a proof of an invariance principle for connected clusters under diffusive scaling.

In two dimensions duality considerations enable a reformulation of these results for supercritical nearest-neighbor random cluster measures, in particular, for nearest-neighbor Potts models in the phase transition regime. Accordingly, we prove that in two dimensions Potts equilibrium crystal shapes are always analytic and strictly convex and that the interfaces between different phases are always diffusive. Thus, no roughening transition is possible in the whole regime where our results apply.

Our results hold under an assumption of exponential decay of finite volume wired connectivities [assumption (1.2) below] in rectangular domains that is conjectured to hold in the whole subcritical regime; the latter is known to be true, in any dimensions, when $q = 1$, $q = 2$, and when $q$ is sufficiently large. In two dimensions assumption (1.2) holds whenever there is an exponential decay of connectivities in the infinite volume measure. By duality, this includes all supercritical nearest-neighbor Potts models with positive surface tension between ordered phases.

Received October 2006; revised July 2007.

[1]Supported in part by Japan Technion Society Research Fund.

[2]Supported in part by the ANR project JC05_42461 POLINTBIO.

*AMS 2000 subject classifications.* Primary 60F15, 60K15, 60K17, 60K35, 82B20, 82B41; secondary 37C30.

*Key words and phrases.* Fortuin–Kasteleyn model, random cluster model, dependent percolation, Ornstein–Zernike behavior for connectivities, invariance principle, Ruelle operator, renormalization, local limit theorems, interfaces in Potts models, strict convexity of Wulff shape, absence of roughening.

**Contents**

**1. Introduction and results.** We recall the definition of random cluster measures $\mathbb{P}$ in Section 1.1 below. Let $x \in \mathbb{Z}^d$ be a distant point. The event $\{0 \leftrightarrow x\}$ means that $x$ is connected to 0 or, equivalently, that the common cluster $\mathbf{C}_{0,x}$ of the origin 0 and $x$ is well defined. In the subcritical situation the probability $\mathbb{P}(0 \leftrightarrow x)$ decays exponentially with $|x|$. The direction-dependent decay exponent $\xi$ is called inverse correlation length. Our main structural result in this paper is a probabilistic description of connected clusters $\mathbf{C}_{0,x}$ in terms of one-dimensional "sausages" of objects, which we call irreducible clusters; see (1.4) below. Under the conditional measure $\mathbb{P}(\cdot \mid 0 \leftrightarrow x)$, the sizes of these irreducible clusters have exponentially decaying tails. This is the original Ornstein–Zernike picture. It was rigorously justified for several particular models either in a perturbative regime in, for example, [5, 19, 20] or for on-axis directions in, for example, [7, 12]. In the latter works the arguments heavily relied on particular microscopic structure of irreducible clusters and on independence between different irreducible components. The main contribution of [8] is a derivation of robust coarse-graining procedures which fits in with essentially any microscopic notion of irreducibility one prefers to choose. The issue of dependence between irreducible components has been resolved for all subcritical finite-range Ising models in [9], where the corresponding statistics have been studied in terms of one-dimensional thermodynamics of Ruelle's type. This is a generic situation. Below we extend the theory of [8, 9] to a large class of finite-range subcritical random cluster measures. The theory applies in any dimension. For two-dimensional nearest-neighbor models duality arguments enable an extension for supercritical interfaces. In particular, our theory leads to a comprehensive fluctuation analysis of interfaces for any two-dimensional nearest-neighbor Potts model in the phase segregation regime.

1.1. *Random-cluster measures.* Let $\mathbf{J} = (J_x)_{x\in\mathbb{Z}^d}$ be a collection of nonnegative real numbers satisfying $J_x = J_{-x}$; we suppose that there exists $R > 0$ such that $J_x = 0$ as soon as $|x| > R$ ($|x|$ denotes Euclidean norm), and that[3] $J_x > 0$ when $|x| = 1$. Let us denote by $\mathcal{E}_\mathbf{J}$ the set of all unordered pairs $(x, y)$ of sites of $\mathbb{Z}^d$ such that $J_{y-x} > 0$; we call bonds the elements of $\mathcal{E}_\mathbf{J}$, and endpoints the two elements of a bond.

Let $\Omega = \{0,1\}^{\mathcal{E}_\mathbf{J}}$ be the set of all bond configurations, and let $\mathcal{F}$ denote the corresponding product $\sigma$-algebra; we say that a bond $e$ is open, respectively closed, in the configuration $\omega \in \Omega$ if $\omega(e) = 1$, respectively $\omega(e) = 0$. A configuration of $\omega \in \Omega$ is identified with the graph with vertex-set $\mathbb{Z}^d$ and edge-set $\{e \in \mathcal{E}_\mathbf{J} : \omega(e) = 1\}$. The maximal connected components of this graph are called open clusters of the configuration (in particular, each isolated site is considered to be an open cluster). We write $x \leftrightarrow y$ for the event that $x$ and $y$ belong to the same open cluster, $A \leftrightarrow B$ if there exist $x \in A$ and $y \in B$ such that $x \leftrightarrow y$, and $x \leftrightarrow \infty$ for the event that the open cluster containing $x$ has infinite cardinality.

[3] This second assumption is purely a matter of convenience; the proofs still apply, modulo straightforward modifications, if this condition is removed, provided that the graph $(\mathbb{Z}^d, \mathcal{E}_\mathbf{J})$ remains connected.

If $\Lambda$ is a finite subset of $\mathbb{Z}^d$, $\Lambda \Subset \mathbb{Z}^d$, we denote by $\mathcal{E}_{\mathbf{J}}^{\Lambda}$ the set of all bonds with both endpoints in $\Lambda$, and by $\mathcal{T}_{\Lambda}$ the $\sigma$-algebra generated by the sets $\{w(e) : e \in \mathcal{E}_{\mathbf{J}} \setminus \mathcal{E}_{\mathbf{J}}^{\Lambda}\}$.

The random cluster measures on $\mathbb{Z}^d$, with coupling constants $\mathbf{J}$ and parameters $\beta \geq 0$ and $q \geq 1$, are all the probability measures $\mathbb{P}$ on $(\Omega, \mathcal{F})$ satisfying the following DLR equations: $\forall A \in \mathcal{F}$ and $\Lambda \Subset \mathbb{Z}^d$,

$$\mathbb{P}(A \mid \mathcal{T}_{\Lambda})(\eta) = \mathbb{P}^{\eta}_{\Lambda;\mathbf{J},\beta,q}(A), \qquad \text{for } \mathbb{P}\text{-a.e. } \eta \in \Omega,$$

where $\mathbb{P}^{\eta}_{\Lambda;\mathbf{J},\beta,q}$ is the probability measure on $(\Omega, \mathcal{F})$ given by

$$\mathbb{P}^{\eta}_{\Lambda;\mathbf{J},\beta,q}(\omega) \stackrel{\text{def}}{=} \begin{cases} (Z^{\eta}_{\Lambda;\mathbf{J},\beta,q})^{-1} \displaystyle\prod_{e \in \mathcal{E}_{\mathbf{J}}^{\Lambda}} p_e^{\omega(e)} (1-p_e)^{1-\omega(e)} q^{\mathfrak{N}_{\Lambda}(\omega)}, & \text{if } \omega \in \Omega^{\eta}_{\Lambda}, \\ 0, & \text{otherwise}, \end{cases}$$

with $p_e \stackrel{\text{def}}{=} 1 - e^{-2\beta J_e}$, $\Omega^{\eta}_{\Lambda} = \{\omega \in \Omega : \omega(e) = \eta(e), \forall e \in \mathcal{E}_{\mathbf{J}} \setminus \mathcal{E}_{\mathbf{J}}^{\Lambda}\}$, and $\mathfrak{N}_{\Lambda}(\omega)$ the number of open clusters of $\omega$ intersecting $\Lambda$.

We write $\mathbb{P}^{\mathrm{w}}_{\Lambda;\mathbf{J},\beta,q}$ when $\eta \equiv 1$ (wired boundary conditions), and $\mathbb{P}^{\mathrm{f}}_{\Lambda;\mathbf{J},\beta,q}$ when $\eta \equiv 0$ (free boundary conditions), and shall generally omit $\mathbf{J}, \beta, q$ from the notation.

The measures $\mathbb{P}^{\mathrm{f}}_{\Lambda}$ and $\mathbb{P}^{\mathrm{w}}_{\Lambda}$ are defined for each finite $\Lambda \subset \mathbb{Z}^d$ and can be viewed as distributions on the set $\{0,1\}^{\mathcal{E}_{\mathbf{J}}^{\Lambda}}$. In a completely similar fashion one constructs measures $\mathbb{P}^{\mathrm{f}}_{\gamma}$ and $\mathbb{P}^{\mathrm{w}}_{\gamma}$ on $\{0,1\}^{\gamma}$ for any finite subgraph (set of bonds) $\gamma$ of $(\mathbb{Z}^d, \mathcal{E}_{\mathbf{J}})$.

A property of random cluster measures that is repeatedly used in the present work is the FKG property: Let $f$ and $g$ be two increasing functions w.r.t. the natural partial order on $\Omega$; then $\mathbb{P}(fg) \geq \mathbb{P}(f)\mathbb{P}(g)$ for any random cluster measure $\mathbb{P}$.

In particular, it follows that, for any increasing sequence $\Lambda \nearrow \mathbb{Z}^d$, the two sequences of measures $\mathbb{P}^{\mathrm{w}}_{\Lambda}$ and $\mathbb{P}^{\mathrm{f}}_{\Lambda}$ converge to random cluster measures $\mathbb{P}^{\mathrm{w}}$ and $\mathbb{P}^{\mathrm{f}}$. Moreover, the following stochastic ordering of measures obtains, $\mathbb{P}^{\mathrm{f}} \preccurlyeq \mathbb{P} \preccurlyeq \mathbb{P}^{\mathrm{w}}$ for any random cluster measure $\mathbb{P}$.

It can also be shown that $\mathbb{P}^{\mathrm{w}}_{\mathbf{J},\beta,q} = \mathbb{P}^{\mathrm{f}}_{\mathbf{J},\beta,q}$ (and therefore, there is a unique random cluster measure), for all but at most countably many values of $\beta$ [15].

In particular, the following definition,

$$\beta_{\mathrm{c}}^1 \stackrel{\text{def}}{=} \inf\{\beta : \mathbb{P}^{*}_{\mathbf{J},\beta,q}(0 \leftrightarrow \infty) > 0\},$$

does not depend on the particular choice $* = \mathrm{w}$ or $* = \mathrm{f}$ (since the event $0 \leftrightarrow \infty$ is increasing). The set $\{\beta : \beta < \beta_{\mathrm{c}}^1\}$ is called a subcritical regime.

It can be shown that there is a unique random cluster measure for each $\beta < \beta_c^1$ [15].

From now on, we suppose that $\beta < \beta_c^1$, and denote the corresponding (unique) random cluster measure by $\mathbb{P}_{\mathbf{J},\beta,q}$, or simply $\mathbb{P}$ when no confusion is possible.

1.1.1. *Inverse correlation length.* The central quantity in our investigation is the connectivity function $\mathbb{P}(0 \leftrightarrow x)$. Associated to the latter is the inverse correlation length: for $x \in \mathbb{R}^d$,

$$\xi(x) \stackrel{\text{def}}{=} - \lim_{k\to\infty} \frac{1}{k} \log \mathbb{P}(0 \leftrightarrow [kx]),$$

where $[y] \in \mathbb{Z}^d$ is the component-wise integer part of $y \in \mathbb{R}^d$.

In the subcritical regime, the connectivity function decays with $|x|$. It is natural to ask whether this decay is actually exponential. We thus introduce another notion of critical point:

$$\beta_c^2 \stackrel{\text{def}}{=} \sup\Big\{\beta : \min_{\vec{n}\in\mathbb{S}^{d-1}} \xi(\vec{n}) > 0\Big\}.$$

Obviously, $\beta_c^2 \le \beta_c^1$, but it is actually believed that these two critical points always coincide.

Conjecture 1. *For any $d \ge 1$, $q \ge 1$, and $\mathbf{J}$ as above,*

$$\beta_c^1 = \beta_c^2.$$

This result is known to hold when $q = 1$ [1], $q = 2$ [2] or $q$ is large enough [17] (in any dimensions).

Let us assume now that $\beta < \beta_c^2$. The function $\xi$ is clearly positively homogeneous; moreover, FKG inequalities imply that $\xi$ is a convex function on $\mathbb{R}^d$, and $\xi$ is obviously finite. It is therefore an equivalent norm on $\mathbb{R}^d$. Two convex bodies play an important role in the sequel: the equi-decay set (or, equivalently, Frank diagram) $\mathbf{U}_\xi$ and the Wulff shape[4] $\mathbf{K}_\xi$, which are defined by[5]

$$\mathbf{U}_\xi \stackrel{\text{def}}{=} \{x \in \mathbb{R}^d : \xi(x) \le 1\}, \qquad \mathbf{K}_\xi \stackrel{\text{def}}{=} \bigcap_{\vec{n}\in\mathbb{S}^{d-1}} \{t \in \mathbb{R}^d : (t,\vec{n})_d \le \xi(\vec{n})\}. \tag{1.1}$$

Notice that these two sets are polar,

$$\mathbf{U}_\xi = \Big\{x \in \mathbb{R}^d : \max_{t\in\mathbf{K}_\xi} (t,x)_d \le 1\Big\}, \qquad \mathbf{K}_\xi = \Big\{t \in \mathbb{R}^d : \max_{x\in\mathbf{U}_\xi} (t,x)_d \le 1\Big\}.$$

[4] This terminology will find its justification when stating our results in the 2D Potts model.

[5] $(x,y)_d$ denotes the scalar product of $x, y \in \mathbb{R}^d$.

Pairs $(x,t) \in \mathbb{R}^d \times \partial \mathbf{K}_\xi$ are called dual if $(t,x)_d = \xi(x)$. The set $\mathbf{K}_\xi$ (or, equivalently, $\mathbf{U}_\xi$) encodes geometrically all the information on the function $\xi$ (in particular, the latter can easily be reconstructed from the set).

1.1.2. *Duality in dimension* 2. We restrict here our attention to the two-dimensional nearest-neighbor case: $J_x = \mathbf{1}_{\{|x|=1\}}$; in this case we simply write $\mathcal{E}$ instead of $\mathcal{E}_{\mathbf{J}}$.

A specific feature will allow us to reinterpret some of our results as pertaining to the $q$-state Potts model in the phase coexistence regime: duality.

Let $\mathbb{Z}^{2,*} = \{(x^*, y^*) : x^* - \frac{1}{2}, y^* - \frac{1}{2} \in \mathbb{Z}\}$ be the dual lattice, and denote by $\mathcal{E}^*$ the nearest-neighbor bonds of $\mathbb{Z}^{2,*}$ (the dual bonds). Let $\Lambda = \{-n, \dots, n\}^2$ and $\Lambda^* = \{-n - \frac{1}{2}, \dots, n + \frac{1}{2}\}^2 \subset \mathbb{Z}^{2,*}$. Random cluster measures on $\Omega^* \stackrel{\text{def}}{=} \{0,1\}^{\mathcal{E}^*}$ are defined exactly as before.

There is a natural mapping from $\Omega \to \Omega^*$: given a configuration $\omega \in \Omega$, the configuration $\omega^* \in \Omega^*$ is obtained by setting, for each dual bond $e^* \in \mathcal{E}^*$, $\omega(e^*) = 1 - \omega(e)$, where $e \in \mathcal{E}$ is the bond dual to $e^*$ (i.e., seen as unit-length line segments in $\mathbb{R}^2$, $e$ is the bond intersecting $e^*$).

Let $p \stackrel{\text{def}}{=} 1 - e^{-2\beta}$, and define $p^*$ as the solution to the equation $p^*/(1 - p^*) = q(1-p)/p$; let $\beta^*$ be such that $p^* = 1 - e^{-2\beta^*}$. Then the following duality relation holds [13]:

$$\mathbb{P}^{\mathrm{f}}_{\Lambda;\beta,q}(\omega) = \mathbb{P}^{\mathrm{w}}_{\Lambda^*;\beta^*,q}(\omega^*).$$

Of course, letting $n \to \infty$, one then also obtains the corresponding statement for infinite-volume measures: $\mathbb{P}^{\mathrm{f}}_{\beta,q}(A) = \mathbb{P}^{\mathrm{w}}_{\beta^*,q}(A^*)$, for all cylindrical events $A$ and $A^* \stackrel{\text{def}}{=} \{\omega^* : \omega \in A\}$.

This duality relates super- and subcritical models. In particular, this will allow us to derive, from our results obtained in the subcritical regime, results about the supercritical regime. Thus, the inverse correlation length $\xi$ of the subcritical model coincides with the surface tension (see, e.g., [17, 18]) $\tau$ of the dual supercritical model. Much more than that: since connected clusters of the subcritical model represent interfaces in the supercritical dual model, our results lead to a comprehensive fluctuation theory of microscopic interfaces which applies to *all* nearest-neighbor two-dimensional random cluster measures (at $q \geq 1$), in particular, it applies to all nearest-neighbor two-dimensional Potts models in the phase segregation regime.

1.1.3. *The basic assumption.* Although the natural region of validity of our results should be the whole region $\{\beta : \beta < \beta_{\mathrm{c}}^2\}$, where $\xi > 0$, it will be necessary to have an *a priori* stronger property. We therefore introduce yet another critical point $\hat{\beta}_{\mathrm{c}}$, defined as the supremum over all values of $\beta$ such that the following holds: Let $\Lambda_N = \{-N, \dots, N\}^d$; then there exist $\nu_0$ and $\nu_1 > 0$ such that, for any $N$,

$$\mathbb{P}^{\mathrm{w}}_{\Lambda_N;\mathbf{J},\beta,q}(0 \leftrightarrow \mathbb{Z}^d \setminus \Lambda_N) \leq \nu_0 e^{-\nu_1 N}. \tag{1.2}$$

Obviously, $\hat{\beta}_{\mathrm{c}} \leq \beta^2_{\mathrm{c}}$, and it is not difficult to show that $\hat{\beta}_{\mathrm{c}} > 0$, for all $d \geq 1$, $q \geq 1$, and $\mathbf{J}$ as above. Again, the following is expected:

Conjecture 2. *For any $d \geq 1$, $q \geq 1$, and $\mathbf{J}$ as above,*

$$\hat{\beta}_{\mathrm{c}} = \beta^2_{\mathrm{c}}.$$

This result has been proved in all the cases mentioned above, where it is known that $\beta^1_{\mathrm{c}} = \beta^2_{\mathrm{c}}$: $q = 1$, $q = 2$ [6], and $q$ sufficiently large [17], in any dimension. Moreover, Conjecture 2 has been verified for any $q \geq 1$ in two dimensions [3]. In the Appendix we give a simple proof of the latter assertion which was kindly explained to us by Reda Messikh.

Our basic assumption throughout this paper is that $\beta < \hat{\beta}_{\mathrm{c}}$.

1.2. *Effective one-dimensional structure of long connected clusters.* From the technical point of view, we encode long connected clusters $\mathbf{C}_{0,x}$ as one-dimensional strings of letters from a countable alphabet of irreducible clusters. The main effort is to verify that the resulting one-dimensional structure falls in the framework of Ruelle's thermodynamic formalism for full shifts with Hölder continuous potentials. Once such identification is accomplished, statistical properties of various local observables comply with classical limit laws. In fact, the only local observable we consider here is the displacement over an irreducible cluster. This leads to an effective random walk representation of $\mathbf{C}_{0,x}$ and provides adequate tools for proving all the main results of the paper. However, we would like to point out that our method provides much more flexibility, and that there are many other local and quasi-local observables of interest.

Thus, before stating precise results, let us give a somewhat informal description of the above mentioned probabilistic structure of clusters $\mathbf{C}_{0,x}$ in terms of a string of irreducible clusters, which lies in the heart of our approach: With each $t \in \partial\mathbf{K}_\xi$, we associate three families $\mathcal{I}^{\mathsf{b}}_t, \mathcal{I}_t$ and $\mathcal{I}^{\mathsf{f}}_t$ of respectively backward irreducible, irreducible and forward irreducible clusters. These families are precisely defined in Section 2.10; see also Figure 1. An important geometric feature is the cone confinement property, which we proceed to describe.

1.2.1. *Geometry of irreducible clusters.* For each $t \in \partial\mathbf{K}_\xi$, define the cone

$$Y = \{y \in \mathbb{Z}^d : (t, y)_d > \kappa_2 \xi(y)\}.$$

Then $t$-irreducible clusters satisfy the following $Y$-cone confinement condition (see Figure 1): For any $\gamma^{\mathsf{b}} \in \mathcal{I}^{\mathsf{b}}_t$, there exists $\mathsf{b} \in \gamma^{\mathsf{b}}$, for any $\gamma^{\mathsf{f}} \in \mathcal{I}^{\mathsf{f}}_t$, there exists $\mathsf{f} \in \gamma^{\mathsf{f}}$ and, finally, for any $\gamma \in \mathcal{I}_t$, there exist $\mathsf{b}, \mathsf{f} \in \gamma$ such that

$$\gamma^{\mathsf{b}} \subset \mathsf{b} - Y, \gamma^{\mathsf{f}} \subset \mathsf{f} + Y \quad \text{and} \quad \gamma \subset (\mathsf{f} + Y) \cap (\mathsf{b} - Y). \tag{1.3}$$

Notice that two irreducible clusters $\gamma$ and $\gamma'$ with $\mathsf{b}(\gamma) = \mathsf{f}(\gamma')$ could be patched together. We use $\gamma \amalg \gamma'$ for the corresponding concatenation.

In (1.4), (1.6) and (1.8) below we state some of the main properties which hold uniformly in $x \in \mathbb{Z}^d$ and $t \in \partial \mathbf{K}_\xi$. Namely, there exist two constants $\kappa_1 = \kappa_1(\mathbf{J}, \beta, q) > 0$ and $\kappa_2 = \kappa_2(\mathbf{J}, \beta, q) \in (0, 1)$ such that:

1.2.2. *Decomposition of* $\mathbf{C}_{0,x}$. Up to probabilities of the order

$$\exp\{-(t, x)_d - \kappa_1 |x|\}$$

clusters $\mathbf{C}_{0,x}$ can be represented as a concatenation,

$$\mathbf{C}_{0,x} = \gamma^{\mathsf{b}} \amalg \gamma_1 \amalg \cdots \amalg \gamma_{\mathcal{N}} \amalg \gamma^{\mathsf{f}}, \tag{1.4}$$

with $\mathcal{N} \geq 1$, $\gamma^{\mathsf{b}} \in \mathcal{I}_t^{\mathsf{b}}$, $\gamma^{\mathsf{f}} \in \mathcal{I}_t^{\mathsf{f}}$ and $\gamma_1, \ldots, \gamma_{\mathcal{N}} \in \mathcal{I}_t$.

Of course, such a statement is meaningful only for $x$'s satisfying $\xi(x) < (t, x)_d + \kappa_1 |x|$.

In Section 3.1 we construct percolation events $\Gamma^{\mathsf{b}}$, $\Gamma_1, \ldots, \Gamma_{\mathcal{N}}$ and $\Gamma^{\mathsf{f}}$, so that the probability of each fixed realization of (1.4) is precisely

$$\mathbb{P}(\Gamma^{\mathsf{b}} \Gamma_1 \cdots \Gamma_{\mathcal{N}} \Gamma^{\mathsf{f}}) = \mathbb{P}(\Gamma^{\mathsf{b}}) \mathbb{P}(\Gamma^{\mathsf{f}}) \frac{\mathbb{P}(\Gamma^{\mathsf{b}} \Gamma_1 \cdots \Gamma_{\mathcal{N}} \Gamma^{\mathsf{f}})}{\mathbb{P}(\Gamma^{\mathsf{b}}) \mathbb{P}(\Gamma^{\mathsf{f}})}. \tag{1.5}$$

1.2.3. *Weights of boundary conditions.* Boundary clusters $\gamma^{\mathsf{b}}$ and $\gamma^{\mathsf{f}}$ cannot be too large: Uniformly in $y \in Y$,

$$\sum_{\gamma^{\mathsf{b}} \ni -y}^{*} \mathbb{P}(\Gamma^{\mathsf{b}}) \leq e^{-(t,y)_d - \kappa_1 |y|} \quad \text{and} \quad \sum_{\gamma^{\mathsf{f}} \ni y}^{*} \mathbb{P}(\Gamma^{\mathsf{f}}) \leq e^{-(t,y)_d - \kappa_1 |y|}, \tag{1.6}$$

where the summation is restricted to backward and forward irreducible clusters satisfying $\mathsf{b}(\gamma^{\mathsf{b}}) = 0$ and, accordingly, $\mathsf{f}(\gamma^{\mathsf{f}}) = 0$.

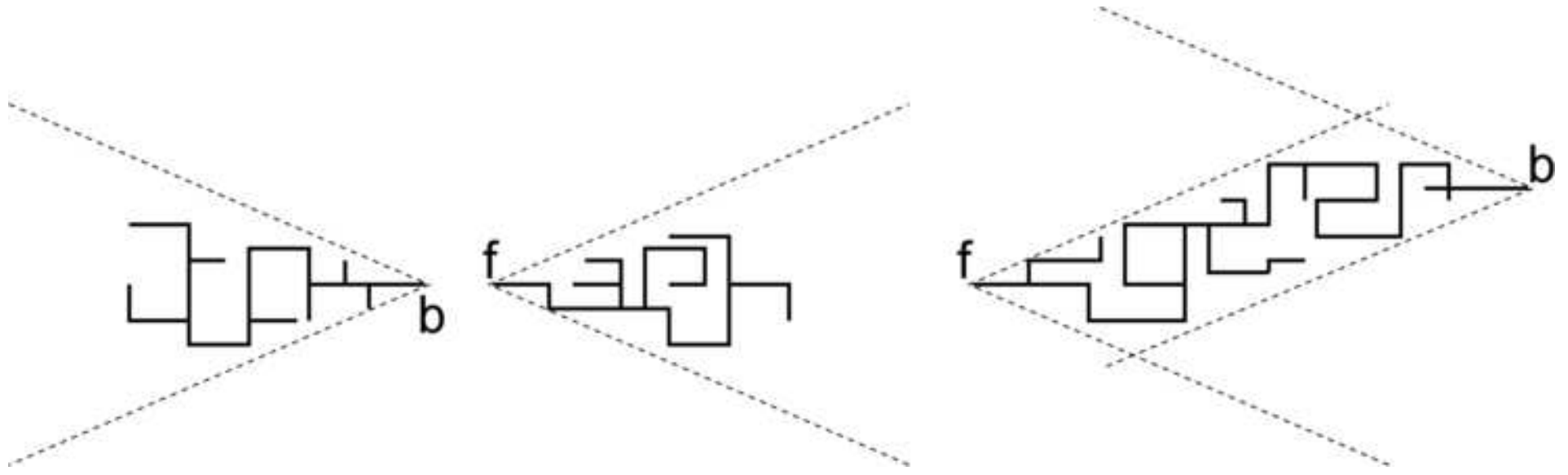

FIG. 1. *Backward ($\mathcal{I}_t^{\mathsf{b}}$) irreducible, forward ($\mathcal{I}_t^{\mathsf{f}}$) irreducible and irreducible ($\mathcal{I}_t$) clusters.*

1.2.4. *Statistics of the effective random walk.* For the purpose of this paper, it is convenient to think about the decomposition (1.4) in terms of an effective random walk with boundary conditions $\gamma^{\mathsf{b}}$ and $\gamma^{\mathsf{f}}$ and steps

$$V(\gamma_1),\dots,V(\gamma_{\mathcal{N}}),$$

where $V(\gamma)=\mathsf{b}-\mathsf{f}$ is the displacement along the irreducible cluster $\gamma\in\mathcal{I}_t$; see Figure 2. In this notation,

$$V(\underline{\gamma})\stackrel{\text{def}}{=}V(\gamma_1)+\cdots+V(\gamma_{\mathcal{N}})$$

is the total displacement along the $\mathcal{N}$-step effective random walk.

For each couple $\gamma^{\mathsf{b}}$ and $\gamma^{\mathsf{f}}$ of boundary conditions, there exist two functions $g(\cdot\mid\gamma^{\mathsf{b}},\gamma^{\mathsf{f}})$ and $\psi_t(\cdot\mid\gamma^{\mathsf{b}},\gamma^{\mathsf{f}})$ acting on strings of irreducible clusters from $\mathcal{I}_t$, such that the fraction on the right-hand side of (1.5) can be represented in the following way: Let $\gamma_1,\dots,\gamma_{\mathcal{N}}$ be irreducible clusters and let $\Gamma_1,\dots,\Gamma_{\mathcal{N}}$ be the corresponding percolation events:

$$\begin{aligned}&\frac{\mathbb{P}(\Gamma^{\mathsf{b}}\Gamma_1\cdots\Gamma_{\mathcal{N}}\Gamma^{\mathsf{f}})}{\mathbb{P}(\Gamma^{\mathsf{b}})\mathbb{P}(\Gamma^{\mathsf{f}})}e^{(t,V(\underline{\gamma}))_d}\\ &\qquad=\exp\{\psi_t(\gamma_1,\dots,\gamma_{\mathcal{N}}\mid\gamma^{\mathsf{b}},\gamma^{\mathsf{f}})+\cdots+\psi_t(\gamma_{\mathcal{N}}\mid\gamma^{\mathsf{b}},\gamma^{\mathsf{f}})\}\\ &\qquad\quad\times g(\gamma_1,\dots,\gamma_{\mathcal{N}}\mid\gamma^{\mathsf{b}},\gamma^{\mathsf{f}}).\end{aligned}\tag{1.7}$$

The potentials $\psi_t$ satisfy

$$\sum_{\gamma_1\in\mathcal{I}_t}e^{\psi_t(\gamma_1,\dots,\gamma_{\mathcal{N}}\mid\gamma^{\mathsf{b}},\gamma^{\mathsf{f}})+\kappa_1|V(\gamma_1)|}<\infty,\tag{1.8}$$

uniformly in $t\in\partial\mathbf{K}_\xi$, in $\gamma^{\mathsf{b}}\in\mathcal{I}_T^{\mathsf{b}}$, $\gamma^{\mathsf{f}}\in\mathcal{I}_T^{\mathsf{f}}$, $\mathcal{N}\in\mathbb{N}$ and in $\gamma_2,\dots,\gamma_{\mathcal{N}}\in\mathcal{I}_t$. In the Bernoulli percolation case ($q=1$) [8] both $\psi_t$ and $g$ depend only on the first element of the string. This corresponds to effective random walks with independent increments which, by (1.8), satisfy Cramér's condition. Accordingly, the fluctuation analysis of connected clusters boils down to an application of classical local limit results of probability theory. For general random cluster models, both the potentials $\psi_t$ and the functions $g$ depend on *all* the variables involved. Nevertheless, in Sections 3.3 and 3.4 we show that they possess appropriate regularity properties, which enable an equally classical local limit analysis along the lines of one-dimensional thermodynamics of full shifts over countable alphabets.

1.3. *Main results.* The geometric description (1.4) of $\mathbf{C}_{0,x}$ in terms of one-dimensional chains of irreducible clusters and the regularity properties of potentials in the induced effective random walk representation (1.7) (see Sections 2 and 3 for precise statements) provide a unified framework for all the results below.

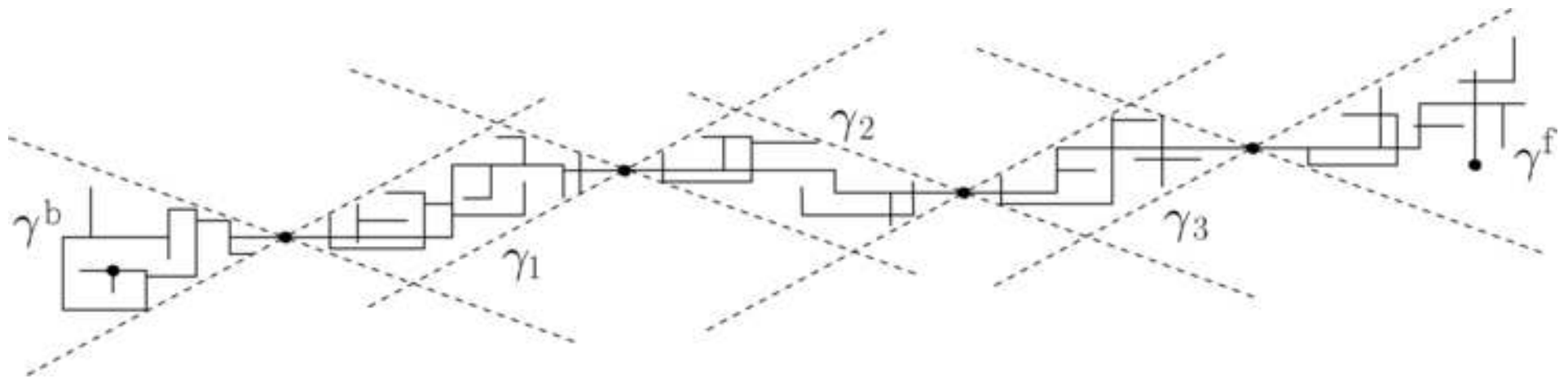

FIG. 2. *Effective random walk representation.*

1.3.1. *Ornstein–Zernike behavior of connectivities.* The principal quantity of interest in this work is the connectivity function,

$$\mathbb{P}(0 \leftrightarrow x),$$

and, more specifically, its asymptotic behavior as $|x| \to \infty$. Our first result, which extends earlier results in the case of self-avoiding walks [16], $q = 1$ [8] and $q = 2$ [9], is that these asymptotics are of the Ornstein–Zernike type.

THEOREM A. *Let $\beta < \hat{\beta}_c$. Then*

$$\mathbb{P}(0 \leftrightarrow x) = \frac{\Psi(\vec{n}_x)}{|x|^{(d-1)/2}} \exp(-\xi(\vec{n}_x)|x|)(1 + o(1)),$$

*uniformly as $|x| \to \infty$. The functions $\Psi$ and $\xi$ are positive, locally analytic functions on $\mathbb{S}^{d-1}$, and $\vec{n}_x \stackrel{\text{def}}{=} x/|x|$.*

As a corollary to the previous theorem, we obtain the following sharp asymptotics for probabilities of exits from sets. For simplicity, we shall formulate such result for cubes $\Lambda_N$ and under the assumption that $\mathbf{J}$ is invariant under all reflection symmetries of $\mathbb{Z}^d$, however, a generalization to generic piece-wise smooth domains $V_N \stackrel{\text{def}}{=} NV$ and general $\mathbf{J}$ is straightforward.[6]

COROLLARY 1.1. *Assume that $\mathbf{J}$ is invariant under all reflections of $\mathbb{Z}^d$. Let $\beta < \hat{\beta}_c$. Then there exists a constant $\psi$, such that*

$$\mathbb{P}(0 \leftrightarrow \mathbb{Z}^d \setminus \Lambda_N) = 2d\psi e^{-N\xi(\mathsf{e}_1)}(1 + o(1)). \tag{1.9}$$

1.3.2. *Geometry of Wulff shapes and equi-decay profiles.* As mentioned above, there is a deep relationship between the geometry of the sets $\mathbf{U}_\xi$ and $\mathbf{K}_\xi$ and the inverse correlation length $\xi$. As a by-product of our analysis, we obtain the following information on the latter.

[6]Here generic means that the minimum of $\xi$ is attained on a discrete set of points of $\partial V$, and $\partial V$ is smooth around these points.

THEOREM B. *Let $\beta < \hat{\beta}_c$. Then $\mathbf{K}_\xi$ has a locally analytic, strictly convex boundary. Moreover, the Gaussian curvature $\chi_\beta$ of $\mathbf{K}_\xi$ is uniformly positive,*

$$\chi_\beta \stackrel{\text{def}}{=} \min_{t\in\partial\mathbf{K}_\xi} \prod_{1}^{d-1} \chi_{\beta,i}(t) > 0,$$

*where $\{\chi_{\beta,i}(t)\}$ are the principal curvatures of $\partial\mathbf{K}_\xi$ at $t$. By duality, $\partial\mathbf{U}_\xi$ is also locally analytic and strictly convex.*

1.3.3. *Invariance principle for long connected clusters.* The local limit result stated in Theorem A can be extended to a full invariance principle for the common cluster $\mathbf{C}_{0,x_n}$ of 0 and $x_n \in \mathbb{Z}^d$ as $x_n \to \infty$. The magnitude of fluctuation depends on the asymptotic direction of $x_n$. Let $x \in \mathbb{S}^{d-1}$ and let $t \in \partial\mathbf{K}_\xi$ be the dual point. Consider a sequence of vertices $x_n = \lfloor nx \rfloor$ and the corresponding sequence of conditional measures $\mathbb{P}_{x,n}(\cdot) \stackrel{\text{def}}{=} \mathbb{P}(\cdot \mid 0 \leftrightarrow x_n)$. By our basic cluster decomposition result, $\mathbf{C}_{0,x_n}$ has the form (1.4) up to $\mathbb{P}_{x,n}$-probabilities of order $e^{-\kappa_1 n}$. Let $\gamma = \{0, u_0, \dots, u_{\mathcal{N}}, x_n\}$ be the trajectory of the corresponding effective random walk, and let $\mathcal{L}_n[\gamma]$ be the linear interpolation through the vertices of $\gamma$. Let $\mathcal{H}_x$ be the $(d-1)$-dimensional hyperplane orthogonal to $x$. Alternatively, $\mathcal{H}_x$ is the tangent space of $\partial\mathbf{K}_\xi$ at $t$. By the cone confinement property (1.3), the intersection number

$$\#(\mathcal{L}_n[\gamma] \cap (r\lfloor nx \rfloor + \mathcal{H}_x)) = 1$$

for every $r \in [0,1]$. Accordingly, there is a natural parametrization of $\mathcal{L}_n[\gamma]$ as a function $\Phi_n : [0,1] \mapsto \mathcal{H}_x$. Namely,

$$\Phi_n(r) = \mathrm{proj}_{\mathcal{H}_x}(\mathcal{L}_n[\gamma] \cap (r\lfloor nx \rfloor + \mathcal{H}_x)).$$

We shall define the diffusive scaling of $\mathbf{C}_{0,x_n}$ in terms of the diffusive scaling of $\Phi_n$. For $r \in [0,1]$, set

$$\phi_n(r) = \frac{1}{\sqrt{n}} \Phi_n(r).$$

Such definition of diffusive scaling of $\mathbf{C}_{0,x_n}$ is justified since by (1.6) and (1.8) the Hausdorff distance between $\mathbf{C}_{0,x_n}$ and $\mathcal{L}_n[\gamma]$ satisfies

$$\lim_{n\to\infty} \mathbb{P}_{x,n}(\mathrm{d_H}(\mathbf{C}_{0,x_n}, \mathcal{L}_n) > K \log n) = 0, \tag{1.10}$$

for some $K = K(\beta)$ large enough.

THEOREM C. *Let $\beta < \hat{\beta}_c$. Then $\{\phi_n(\cdot)\}$ weakly converges under $\{\mathbb{P}_{x,n}\}$ to the distribution of*

$$(\sqrt{\chi_{\beta,1}} B_1(\cdot), \dots, \sqrt{\chi_{\beta,d-1}} B_{d-1}(\cdot)),$$

*where $B_1, \dots, B_{d-1}$ are independent standard Brownian bridges on $[0,1]$.*

1.3.4. *Two dimensional Potts models.* In the special case of two-dimensional models with nearest-neighbor interactions $J_x = \mathbf{1}_{\{|x|=1\}}$, duality allows us to reinterpret the preceding results as results about interfaces in the supercritical random cluster model, as explained in Section 1.1.2. In particular, when $q$ is an integer, this provides a detailed description of interfaces in the $q$-state Potts model in the phase coexistence regime. Let us briefly recall the definition of this model. Let $q \geq 2$ be an integer, and $\Lambda \Subset \mathbb{Z}^2$. The nearest-neighbor $q$-states Potts measure in $\Lambda$ with boundary condition $\bar\sigma \in \{1,\dots,q\}^{\mathbb{Z}^2}$, at inverse temperature $\beta$, is the probability measure on the product $\sigma$-algebra of $\{1,\dots,q\}^{\mathbb{Z}^2}$ given by

$$\Xi_\Lambda^{\bar\sigma}(\sigma) = \begin{cases} (Z_\Lambda^{\bar\sigma})^{-1} \displaystyle\prod_{\{i,j\}\cap\Lambda\neq\varnothing, |i-j|=1} e^{\beta\delta_{\sigma_i,\sigma_j}}, & \text{if } \sigma \equiv \bar\sigma \text{ off } \Lambda, \\ 0, & \text{otherwise.} \end{cases}$$

Since the inverse correlation length of the subcritical model coincides with the surface tension of the supercritical model, the set $\mathbf{K}_\xi$ really corresponds to the Wulff shape associated to the dual model (see [4] for a review on Wulff construction). Therefore, observing that, as already mentioned, $\hat\beta_c$ is known to coincide with $\beta_c^2$, we see that we can restate Theorem B as follows.

THEOREM D. *Potts models in two dimensions do not have roughening transition. In other words, for any $q \in \mathbb{N}$ and at any temperature $\beta$ at which the surface tension $\xi$ of the two-dimensional nearest-neighbor $q$-state Potts model is positive, the corresponding equilibrium crystal (Wulff) shape $\mathbf{K}_\xi$ is strictly convex and has an analytic boundary $\partial\mathbf{K}_\xi$ with strictly positive curvature $\chi_\beta$.*

This result was previously known only in the cases $q = 1$ and $q = 2$.

Let us now consider the 2D $q$-state Potts model in $\Lambda_N = \{-N,\dots,N\}^2$, with $\vec{n}$-Dobrushin boundary condition, $\vec{n} \in \mathbb{S}^1 : \bar\sigma_i = 1$ if $(\vec{n}, i)_d \geq 0$, and $\bar\sigma_i = 2$ if $(\vec{n}, i)_d < 0$. For definiteness, let us assume that $(\vec{n}, \vec{e}_2)_d \geq 1/\sqrt{2}$, and write $\Xi_\Lambda^{\vec{n}}$ for the corresponding measure. Under such a boundary condition, we can define unambiguously the interface, for any configuration $\sigma$ having positive probability under $\Xi_\Lambda^{\vec{n}}$: Among all connected components of all dual edges separating disagreeing spins, only one has infinite cardinality; the intersection of the latter with the box $\Lambda$ defines the interface (see Figure 3). It turns out (see the proof of the following theorem) that the interface as defined here is a subset of an open cluster in the subcritical dual of the random cluster model to which it is associated through the standard coupling [21], and thus can be approximated by the same effective random walk. Denoting as before by $\Phi_N$ its linear interpolation, and by $\phi_N$ the diffusive scaling of the latter, we arrive at the following theorem (notice that this is not an immediate corollary of Theorem C because the corresponding random cluster model is defined in a vertical strip, and one must therefore take care of the boundary effects).

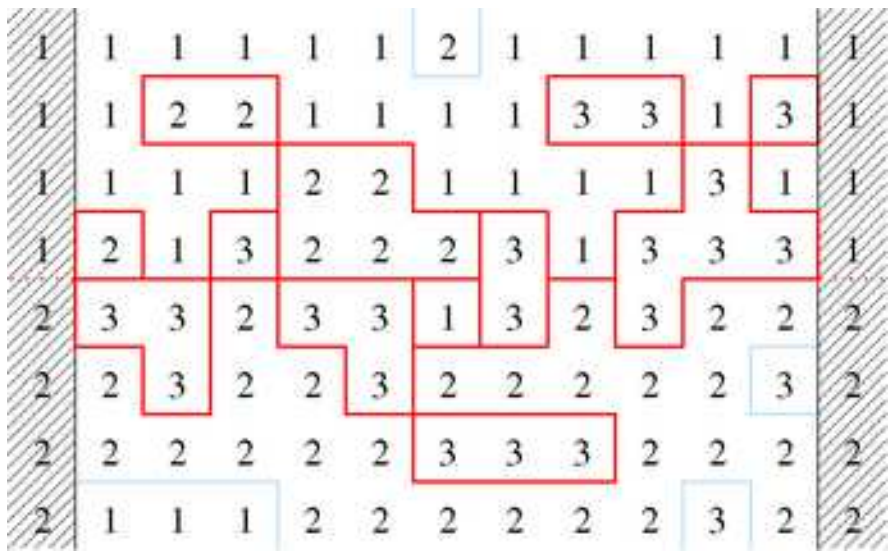

FIG. 3. *The interface of the 2D nearest-neighbor three state Potts model.*

THEOREM E. *Consider the $q$-states nearest-neighbor Potts model with $\vec{n}$-Dobrushin boundary condition in the box $\Lambda_N$, at inverse temperature $\beta$ such that the surface tension is positive. Then $\{\phi_N(\cdot)\}$ weakly converges under $\Xi_\Lambda^{\vec{n}}$ to the distribution of*

$$\{\sqrt{\chi_\beta} B(\cdot)\},$$

*where $B$ is the standard Brownian bridge on $[0,1]$ and $\chi_\beta$ is the curvature of $\partial \mathbf{K}_\xi$ computed at the point $t \in \partial \mathbf{K}_\xi$ dual to $\vec{n}$.*

The case $q = 2$ has already been treated in [14].

Furthermore, the theory we develop also leads to a nonperturbative microscopic justification of the Wulff construction for two-dimensional Potts models on the DKS-level, and also implies that the inverse correlation length is an analytic function of $\beta$. However, we relegate the corresponding discussion to a future paper.

1.4. *Structure of the paper.* As it is indicated in Section 1.2 our main new result in this work is a robust stochastic-geometric description of long connected clusters $\mathbf{C}_{0,x}$, which is valid for all subcritical random cluster measures satisfying assumption (1.2).

In Section 2, a finite scale renormalization analysis of the typical geometry of connected clusters implies that $\mathbf{C}_{0,x}$ has, up to negligible probabilities, an essentially one-dimensional structure, which can be conveniently visualized as a chain of irreducible clusters. The crucial model-oriented input for an adaptation of the general coarse-graining techniques developed in [8, 9] is a mixing estimate on exponential decay of connectivities under various boundary conditions. This estimate is derived in Section 2.1 as a direct consequence of our assumption (1.2). The rest of Section 2 is devoted then to the construction of the irreducible decomposition (1.4) of the cluster $\mathbf{C}_{0,x}$.

The second step of the proof is a refined local limit analysis of chains of irreducible clusters. Section 3 is devoted to a representation of $\mathbb{P}(0 \leftrightarrow x)$ in terms of thermodynamics of Ruelle operators with uniformly Lipschitz continuous potentials for full shifts over countable alphabets (of irreducible clusters). The main representation formula (3.3) is derived in Section 3.1. A mapping to the thermodynamic formalism is explained in Section 3.3. The key result is Theorem 3.1 of Section 3.4 which asserts that assumption (1.2) implies appropriate regularity properties of potentials.

Once Theorem 3.1 is established, the general theory of [9] applies and all the results announced in Section 1.3 follow in a relatively straightforward fashion. In Section 4 we briefly explain the corresponding proofs along the lines of [9, 10, 14].

Finally, in the Appendix we give a simple proof that in two dimensions assumption (1.2) holds whenever the inverse correlation length is positive. (This proof was explained to us by Reda Messikh and we publish it here with his kind permission.)

*Remark on constants.* Positive constants $c_1, c_2, \ldots$ are used in the intermediate computations and are updated with each section. Constants $\kappa_1, \kappa_2$ and $\nu_1, \nu_2, \ldots$ appear in our principal estimates and are fixed throughout the paper. Finally, finite scale renormalization procedures we employ depend on three parameters $K$, $r$ and $\delta$ as it is described in Section 2 below.

## 2. Geometry of typical clusters.

2.1. *Notation and basic decay estimate.* The geometry of the problem is conveniently spelled out in terms of the equi-decay set $\mathbf{U}_\xi$ defined in (1.1). Let us fix a number $r > 0$ and a finite scale $K > 0$, both of which will be later chosen large enough, depending on $\beta$ and on the cone-opening parameter $\delta$ which we shall define in Section 2.6 below. For any $y \in \mathbb{Z}^d$, set

$$\mathbf{B}_K(y) \stackrel{\text{def}}{=} (y + K \cdot \mathbf{U}_\xi) \cap \mathbb{Z}^d \quad \text{and} \quad \bar{\mathbf{B}}_K(y) \stackrel{\text{def}}{=} \mathbf{B}_{K + r \log K}(y).$$

For $y = 0$, we shall use the shorthand notation

$$\mathbf{B}_K \stackrel{\text{def}}{=} \mathbf{B}_K(0) \quad \text{and, accordingly,} \quad \bar{\mathbf{B}}_K \stackrel{\text{def}}{=} \bar{\mathbf{B}}_K(0).$$

Furthermore, for any subset $A \subseteq \mathbb{Z}^d$, set

$$\bar{A}_K \stackrel{\text{def}}{=} \bigcup_{y \in A} \mathbf{B}_{r \log K}(y).$$

We shall also use the notation

$$\partial_R A \stackrel{\text{def}}{=} \{y \in \mathbb{Z}^d \setminus A : d(y, A) \le R\}$$

for the $R$-outer boundary of $A$, where $R$ denotes the range of the process (i.e., the length of the longest possible bond). Given a subset $A \subseteq \mathbb{Z}^d$ and two points $x \in A$ and $y \in A \cup \partial_R A$, we shall use $\{x \stackrel{A}{\longleftrightarrow} y\}$ to denote the event that $x$ and $y$ are connected by an open path $\lambda : x \mapsto y$, such that all the vertices of $\lambda$, with the possible exception of the terminal point $y$ itself, belong to $A$.

There are two parameters $K$ and $r$ (which enter the picture through the definitions of $\mathbf{B}_K, \bar{\mathbf{B}}_K$ and $\bar{A}_K$ above) to play with. Let us say that a certain quantity is of order $o_K(1)$ if for every power $m > 0$ one can find $r > 0$ and a scale $K_0$, such that for every $K \geq K_0$ this quantity is bounded above by $1/K^m$.

The next proposition is the only place where we rely on our basic estimate (1.2) (beyond just assuming exponential decay of connectivities).

Proposition 2.1. *For any subset $A \subseteq \mathbf{B}_K$ and every vertex $y \in \partial_R \mathbf{B}_K$,*

$$\mathbb{P}^{\mathrm{w}}_{\bar{A}_K}(0 \stackrel{A}{\longleftrightarrow} y) \leq e^{-K}(1 + o_K(1)). \tag{2.1}$$

Proof. The proof is rather standard (see Figure 4): By the FKG property of $\mathbb{P}$, we can decompose $\mathbb{P}^{\mathrm{w}}_{\bar{A}_K}(0 \stackrel{A}{\longleftrightarrow} y)$ as follows:

$$\begin{aligned}
\mathbb{P}^{\mathrm{w}}_{\bar{A}_K}(0 \stackrel{A}{\longleftrightarrow} y) &= \mathbb{P}^{\mathrm{w}}_{\bar{A}_K}(0 \stackrel{A}{\longleftrightarrow} y; \partial_R A \leftrightarrow \partial_R \bar{A}_K) \\
&\quad + \mathbb{P}^{\mathrm{w}}_{\bar{A}_K}(0 \stackrel{A}{\longleftrightarrow} y; \partial_R A \not\leftrightarrow \partial_R \bar{A}_K) \\
&\leq \mathbb{P}^{\mathrm{w}}_{\bar{A}_K}(0 \stackrel{A}{\longleftrightarrow} y)\mathbb{P}^{\mathrm{w}}_{\bar{A}_K \setminus A}(\partial_R A \leftrightarrow \partial_R \bar{A}_K) \\
&\quad + \mathbb{P}^{\mathrm{f}}_{\bar{A}_K}(0 \stackrel{A}{\longleftrightarrow} y).
\end{aligned}$$

By assumption (1.2),

$$\mathbb{P}^{\mathrm{w}}_{\bar{A}_K \setminus A}(\partial_R A \leftrightarrow \partial_R \bar{A}_K) = o_K(1).$$

On the other hand,

$$\mathbb{P}^{\mathrm{f}}_{\bar{A}_K}(0 \stackrel{A}{\longleftrightarrow} y) \leq \mathbb{P}(0 \stackrel{A}{\longleftrightarrow} y) \leq e^{-K}. \qquad \square$$

2.2. *Skeletons.* In this subsection we shall develop a coarse-grained picture of the common cluster $\mathbf{C}_{0,x}$ of 0 and $x \in \mathbb{Z}^d$. Fix a scale $K$ and a number $r$ as in the beginning of Section 2.1.

Given a realization $\mathbf{C}$ of $\mathbf{C}_{0,x}$, we are going to construct a tree $\mathcal{T}(\mathbf{C})$, with vertices $x_0 = 0, x_1, \ldots, x_{N(\mathcal{T})}$ using the following procedure.

Step 1. Set $x_0 = 0$, $\mathcal{C} = \bar{\mathbf{B}}_K(x_0)$, $i = 1$.

Step 2. If there is at least one vertex $y \in \partial_R \mathcal{C}$ such that

$$y \stackrel{\mathbf{C} \setminus \mathcal{C}}{\longleftrightarrow} \partial_R \mathbf{B}_K(y) \setminus \mathcal{C},$$

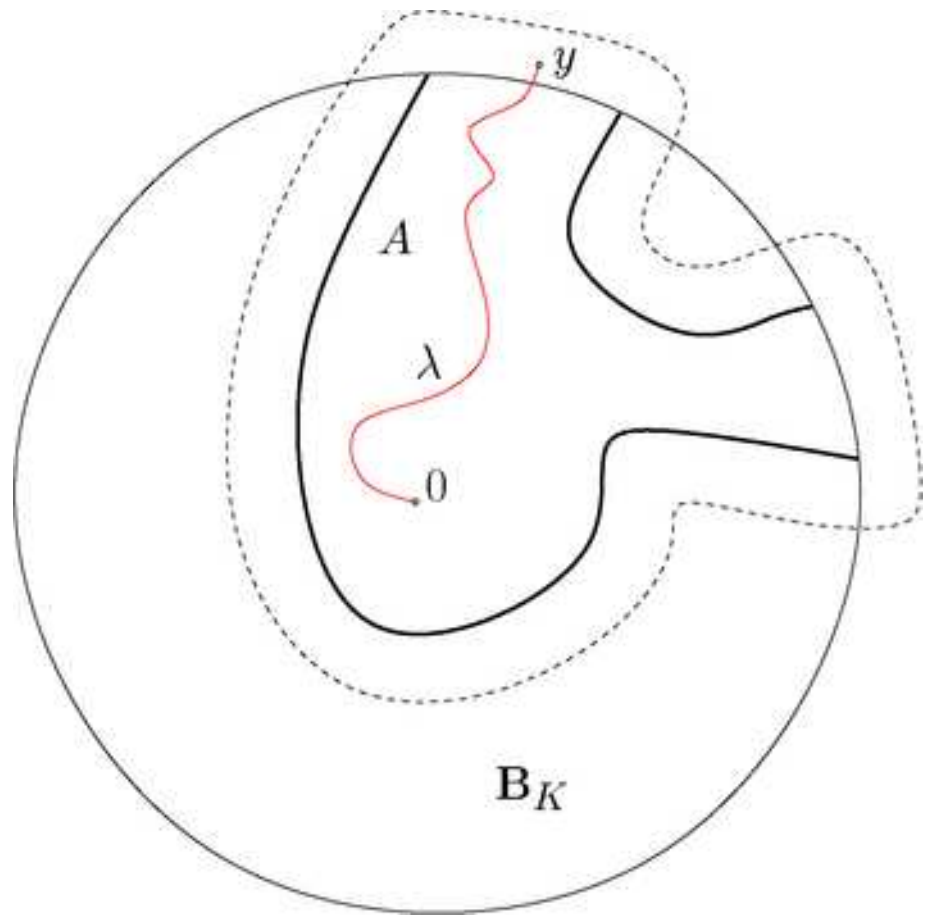

FIG. 4. *The event* $\{0 \overset{A}{\longleftrightarrow} y\}$. *The dashed curve is the boundary* $\partial_R \bar{A}_K$.

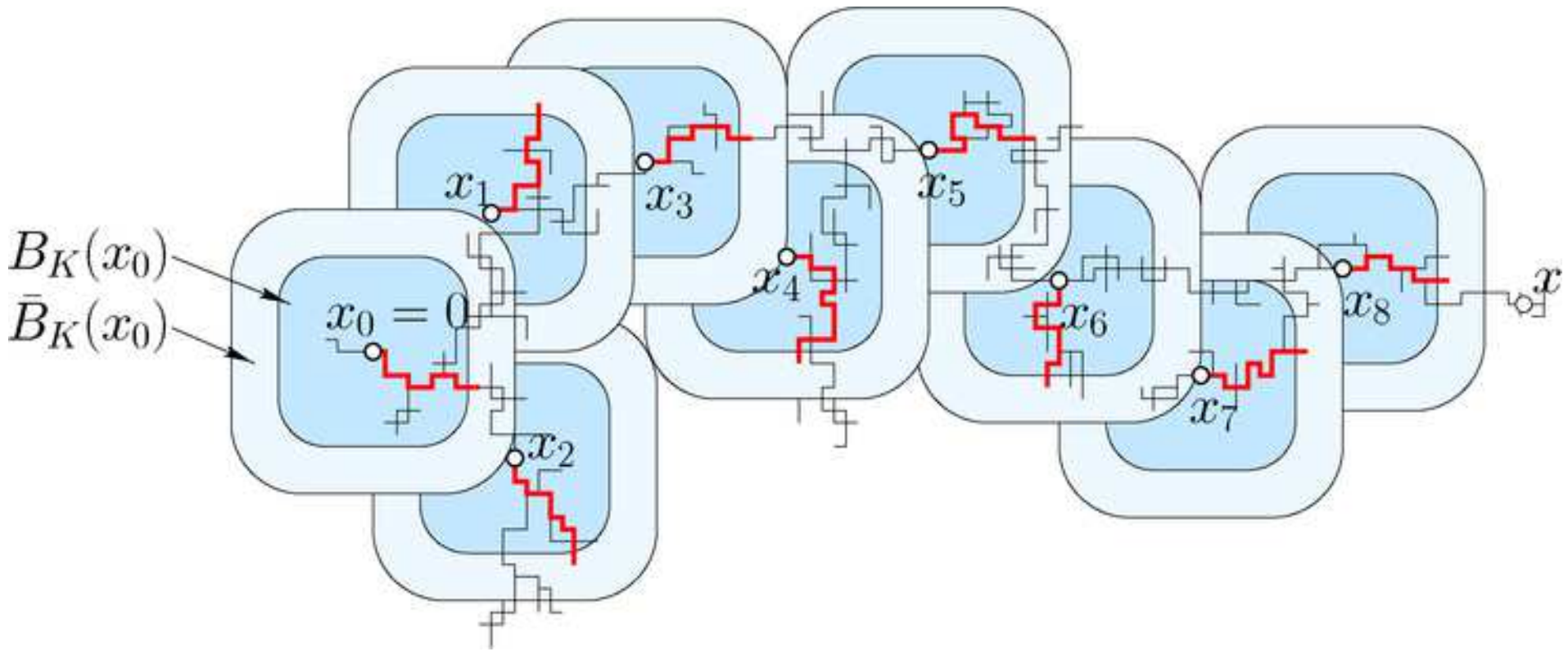

FIG. 5. *The tree skeleton associated to* $\mathbf{C}_{0,x}$.

then choose $y^*$ to be the minimal (in the lexicographic order) such vertex and go to Step 3. Otherwise stop the procedure.

STEP 3. Set $x_i = y^*$. Update $\mathcal{C} \to \mathcal{C} \cup \bar{\mathbf{B}}_K(x_i)$ and $i \to i+1$. Go to Step 2.

This procedure produces the vertices of a tree $\mathcal{T}(\mathbf{C})$. Clearly,

$$\bigcup_i \bar{\mathbf{B}}_K(x_i) \text{ is connected, and } \mathbf{C} \subseteq \bigcup_i \bar{\mathbf{B}}_{2K}(x_i).$$

An example of a tree skeleton of a cluster is depicted on Figure 5.

We now give a rule to construct the edges of the tree: By definition, the vertex $x_i$ is connected to the vertex of

$$\{x_j, j < i : x_i \in \partial_R \bar{\mathbf{B}}_K(x_j)\},$$

which has the smallest index. We define the probability of such a tree $\mathcal{T}$ as

$$\mathbb{P}(\mathcal{T}) \stackrel{\text{def}}{=} \sum_{\mathbf{C} \sim \mathcal{T}} \mathbb{P}(\mathbf{C}_{0,x} = \mathbf{C}),$$

where the sum is over all clusters containing both 0 and $x$ compatible with the above procedure (in the sense that the tree $\mathcal{T}$ can be obtained from the cluster using the procedure).

By Proposition 2.1, we have the following bound on the probability of observing a particular tree $\mathcal{T}$: For each $i = 0, 1, \ldots, N(\mathcal{T})$, set

$$A^i = \mathbf{B}_K(x_i) \setminus \bigcup_{j<i} \bar{\mathbf{B}}_K(x_j).$$

By our construction of $\mathcal{T}$,

$$\begin{aligned} \mathbb{P}(\mathcal{T}) &\le \mathbb{P}(x_i \stackrel{A^i}{\longleftrightarrow} \partial_R \mathbf{B}_K(x_i) \text{ for } i = 0, 1, \ldots, N(\mathcal{T})) \\ &\le \prod_{i=0}^{N(\mathcal{T})} \mathbb{P}^{\mathrm{w}}_{\bar{A}^i_K}(x_i \stackrel{A^i}{\longleftrightarrow} \partial_R \mathbf{B}_K(x_i)) \qquad (2.2) \\ &\le \exp(-K(1 - o_K(1))(N(\mathcal{T}) + 1)), \end{aligned}$$

where we have used (2.1) on the last step.

We are now going to relabel the vertices of a tree $\mathcal{T}$ in the following way: $\mathcal{T}$ can be conveniently split into a trunk $\mathfrak{t} = (\mathfrak{t}_0, \ldots, \mathfrak{t}_{N(\mathfrak{t})})$, and a family of disjoint branches $\mathfrak{B} = (\mathfrak{b}^1, \mathfrak{b}^2, \ldots)$; $\mathcal{T} = (\mathfrak{t}, \mathfrak{B})$. The trunk is defined as the path in $\mathcal{T}$ connecting $\mathfrak{t}_0 = x_0$ to $\mathfrak{t}_N = x_{\mathrm{F}}$, where the latter vertex is defined to be the vertex $x_i$ of $\mathcal{T}$ with the smallest index $i$, having the property that $x \in \bar{\mathbf{B}}_{2K}(x_i)$. Each branch $\mathfrak{b}^l$ is then a connected sub-tree rooted at some vertex of the trunk. The symbols $N(\mathfrak{t})$ and $N(\mathfrak{B}) = \sum_l N(\mathfrak{b}^l)$ are reserved for the total number of vertices of $\mathfrak{t}$ and $\mathfrak{B}$ respectively. Of course, $N(\mathcal{T}) = N(\mathfrak{t}) + N(\mathfrak{B})$ in this notation.

The probability of a trunk $\mathfrak{t}$ is defined as

$$\mathbb{P}(\mathfrak{t}) \stackrel{\text{def}}{=} \sum_{\mathbf{C} \sim \mathfrak{t}} \mathbb{P}(\mathbf{C}_{0,x} = \mathbf{C}),$$

where the compatibility $\mathbf{C} \sim \mathfrak{t}$ means that $\mathfrak{t}$ is a trunk compatible with this construction, and similarly for the leaves. The main idea of the renormalization approach is that on large enough finite scales $K$ typical trees have a simple geometry: Most of the increments along the trunk $\mathfrak{t}$ point into the direction of the target vertex $x$, whereas the branches $\mathfrak{B}$ are very short and sparse. Thus, our estimates on probabilities of trees follow the scheme

$$\mathbb{P}((\mathfrak{t},\mathfrak{B})) = \mathbb{P}(\mathfrak{t})\mathbb{P}(\mathfrak{B} \mid \mathfrak{t}). \tag{2.3}$$

Our bounds on $\mathbb{P}(\mathfrak{B} \mid \mathfrak{t})$ are crude: One pays at least $e^{-K(1-o_K(1))}$ for each additional vertex of $\mathcal{T}$ and, as we shall see in Sections 2.3 and 2.4, at large enough scales $K$ such a price beats the entropy of different tree patterns. The main point is to control the forward geometry of typical trunks $\mathfrak{t}$, which we shall do with the help of the surcharge function introduced in Section 2.5.

In all the estimates below probabilities of various tree patterns will be tested against an a priori bound

$$\mathbb{P}(0 \leftrightarrow x) \asymp e^{-\xi(x)}. \tag{2.4}$$

2.3. *Typical trunks do not contain too many vertices.* The first observation is that the trunk $\mathfrak{t}$ of $\mathbf{C}_{0,x}$ cannot contain too many vertices.

LEMMA 2.1. *There exist $r$, a scale $K_0$, and two constants $c_1$ and $c_2 > 0$ such that*

$$\mathbb{P}\biggl(N(\mathfrak{t}) > c_1\frac{|x|}{K} \biggm| 0 \leftrightarrow x\biggr) \le e^{-c_2|x|},$$

*uniformly in $x$ and $K > K_0$.*

PROOF. The number of trunks of $m$ vertices is bounded above by $(c_3(d,R)K^{d-1})^m$, and therefore, using (2.2),

$$\mathbb{P}\biggl(N(\mathfrak{t}) > c_1\frac{|x|}{K}\biggr) \le \sum_{m > c_1|x|/K} \exp\biggl(-mK\biggl(1 - o_K(1) - \frac{\log(c_3K^{d-1})}{K}\biggr)\biggr),$$

and the conclusion follows by (2.4). □

2.4. *Typical trees do not contain too many branches.*

LEMMA 2.2. *Let $\kappa > 0$. There exist $r$, $K_0(\kappa)$, and $c_4 = c_4(\kappa) > 0$ such that*

$$\mathbb{P}\biggl(N(\mathfrak{B}) > \kappa\frac{|x|}{K} \biggm| 0 \leftrightarrow x\biggr) \le e^{-c_4\kappa|x|},$$

*uniformly in $x$ and $K > K_0$.*

PROOF. By Lemma 2.1, we can assume that $N(\mathfrak{t}) \le c_1\frac{|x|}{K}$. We shall call such trunks *admissible*. Given an admissible trunk, the number of trees with $N(\mathfrak{B}) = b > \kappa\frac{|x|}{K}$ is bounded above by (see Lemma 2.3 in [8])

$$\exp\biggl(c_5 b\biggl(\log\biggl(\frac{1}{\kappa}\biggr) + \log K\biggr)\biggr).$$

Therefore, for any admissible trunk $\mathfrak{t}$ with $N(\mathfrak{t}) \le c_1 \frac{|x|}{K}$, one has

$$\mathbb{P}\Big(N(\mathfrak{B}) > \kappa\frac{|x|}{K}\,\Big|\,\mathfrak{t}\Big) \le \sum_{b>\kappa|x|/K} \exp\Big(-bK\Big(1 - c_5\frac{\log(1/\kappa)+\log K}{K}\Big)\Big),$$
$$\le e^{-\kappa|x|/2},$$

provided $K$ is sufficiently large. □

2.5. *Surcharge function.* On large enough renormalization scales $K$, most of the increments along typical trunks $\mathfrak{t}$ point in the direction of $x$. More precisely, let $t \in \partial\mathbf{K}_\xi$ be a *dual* vector such that $(t,x)_d = \xi(x)$. Then the probability of a trunk $\mathfrak{t} = (\mathfrak{t}_0, \ldots, \mathfrak{t}_N)$ is quantitatively measured using the surcharge functional

$$\mathfrak{s}_t(\mathfrak{t}) \stackrel{\text{def}}{=} \sum_{l=1}^{N} \mathfrak{s}_t(\mathfrak{t}_l - \mathfrak{t}_{l-1}),$$

where the surcharge function $\mathfrak{s}_t : \mathbb{R}^d \mapsto \mathbb{R}_+$ is defined via

$$\mathfrak{s}_t(y) \stackrel{\text{def}}{=} \xi(y) - (t,y)_d.$$

By construction, $\mathfrak{s}_t(x) = 0$.

Now, for all $y \in \partial_R \bar{\mathbf{B}}_K$,

$$K \ge \xi(y) - c_6 r\log K = (t,y)_d + \mathfrak{s}_t(y) - c_6 r \log K.$$

Since by construction $\mathfrak{t}_l \in \partial_R\bar{\mathbf{B}}_K(\mathfrak{t}_{l-1})$ for every $l = 1, \ldots, N(\mathfrak{t})$, and since $x \in \bar{\mathbf{B}}_{2K}(\mathfrak{t}_N)$, we infer, setting $n = N(\mathfrak{t})$, that

$$nK \ge \sum_1^n \xi(\mathfrak{t}_l - \mathfrak{t}_{l-1}) - c_6 rn\log K = (t, \mathfrak{t}_N)_d + \mathfrak{s}_t(\mathfrak{t}) - c_6 rn \log K$$
$$\ge \xi(x) + \mathfrak{s}_t(\mathfrak{t}) - c_7 rn\log K,$$

where on the last step we have tacitly assumed that $n \gg K$. Reasoning as in (2.2), we arrive to the following upper bound:

$$\mathbb{P}(\mathfrak{t}) \le \exp\Big(-\xi(x) - \mathfrak{s}_t(\mathfrak{t}) + c_8\frac{r\log K}{K}|x|\Big), \tag{2.5}$$

for every admissible trunk $\mathfrak{t}$.

We are ready to state our crucial surcharge inequality, which enables a control of the forward geometry of typical trunks of $\mathbf{C}_{0,x}$:

Lemma 2.3. *Let $\varepsilon > 0$. There exist $r$ and $K_0(\varepsilon)$ such that*

$$\mathbb{P}(\mathfrak{s}_t(\mathfrak{t}) > 2\varepsilon|x|, 0 \leftrightarrow x) \le e^{-\xi(x) - \varepsilon|x|},$$

*uniformly in $x \in \mathbb{Z}^d$, vectors $t \in \partial\mathbf{K}_\xi$ dual to $x$, and renormalization scales $K > K_0$.*

PROOF. By Lemma 2.1, we can restrict attention to admissible trunks. The number of the latter is bounded above by

$$\exp\biggl(c_9(d,R)\frac{|x|\log K}{K}\biggr).$$

Therefore, in view of (2.5),

$$\begin{aligned}
&\mathbb{P}(\mathfrak{s}_t(\mathfrak{t}) > 2\varepsilon|x|, N(\mathfrak{t}) \le c_1|x|/K, 0 \leftrightarrow x)\\
&\quad\le \exp\biggl(-|x|\biggl(2\varepsilon - (c_8 r + c_9)\frac{\log K}{K}\biggr) - \xi(x)\biggr),
\end{aligned}$$

and the conclusion follows once $K$ is large enough. □

2.6. *Cone points of trunks.* Let $\delta \in (0,\frac{1}{3})$. For any $t \in \partial\mathbf{K}_\xi$, we define the forward cone by

$$Y_\delta^>(t) \stackrel{\text{def}}{=} \{x \in \mathbb{Z}^d : \mathfrak{s}_t(x) < \delta\xi(x)\},$$

and the backward cone by $Y_\delta^<(t) \stackrel{\text{def}}{=} -Y_\delta^>(t)$.

Given a trunk $\mathfrak{t} = (\mathfrak{t}_0,\dots,\mathfrak{t}_{N(\mathfrak{t})})$, we say that $\mathfrak{t}_k$ is a $(t,\delta)$-forward cone point of $\mathfrak{t}$ if

$$\{\mathfrak{t}_{k+1},\dots,\mathfrak{t}_{N(\mathfrak{t})}\} \subset \mathfrak{t}_k + Y_\delta^>(t),$$

while $\mathfrak{t}_k$ is a $(t,\delta)$-backward cone point of $\mathfrak{t}$ if

$$\{\mathfrak{t}_0,\dots,\mathfrak{t}_{k-1}\} \subset \mathfrak{t}_k + Y_\delta^<(t) = \mathfrak{t}_k - Y_\delta^>(t).$$

$\mathfrak{t}_k$ is then said to be a $(t,\delta)$-cone point of $\mathfrak{t}$ if it is both a $(t,\delta)$-forward cone point and a $(t,\delta)$-backward cone point of $\mathfrak{t}$.

If $\mathfrak{t}$ contains points that are not forward cone points, we define

$$\begin{aligned}
l_1^> &= \min\{j : \mathfrak{t}_j \text{ is not a } (t,\delta)\text{-forward cone point of } \mathfrak{t}\},\\
r_1^> &= \min\{j > l_1 : \mathfrak{t}_j - \mathfrak{t}_{l_1} \notin Y_\delta^>(t)\},\\
l_2^> &= \min\{j \ge r_1 : \mathfrak{t}_j \text{ is not a } (t,\delta)\text{-forward cone point of } \mathfrak{t}\},\\
r_2^> &= \min\{j > l_2 : \mathfrak{t}_j - \mathfrak{t}_{l_2} \notin Y_\delta^>(t)\},\\
&\dots
\end{aligned}$$

and, similarly, if $\mathfrak{t}$ contains points that are not backward cone points,

$$\begin{aligned}
l_1^< &= \max\{j : \mathfrak{t}_j \text{ is not a } (t,\delta)\text{-backward cone point of } \mathfrak{t}\},\\
r_1^< &= \max\{j < l_1 : \mathfrak{t}_j - \mathfrak{t}_{l_1} \notin Y_\delta^<(t)\},\\
l_2^< &= \max\{j \le r_1 : x_j \text{ is not a } (t,\delta)\text{-backward cone point of } \mathfrak{t}\},\\
r_2^< &= \max\{j < l_2 : \mathfrak{t}_j - \mathfrak{t}_{l_2} \notin Y_\delta^<(t)\},\\
&\dots
\end{aligned}$$

We then say that $\mathfrak{t}_j$ is a $(t,\delta)$-marked point of $\mathfrak{t}$ if

$$j \in \bigvee_k \{l_k^>, \ldots, r_k^> - 1\} \cup \bigvee_k \{r_k^< + 1, \ldots, l_k^<\},$$

where $\bigvee$ denotes the disjoint union.

The next lemma, the proof of which is the same as the proof of Lemma 2.2 in [9], controls the number of marked points, $\#_{t,\delta}^{\text{mark}}(\mathfrak{t})$, in terms of the surcharge function.

Lemma 2.4. *Let $\delta \in (0, \frac{1}{3})$ be fixed. The surcharge cost $\mathfrak{s}_t(\mathfrak{t})$ is controlled in terms of the number of marked points as*

$$\mathfrak{s}_t(\mathfrak{t}) \geq c_{10} \delta K \#_{t,\delta}^{\text{mark}}(\mathfrak{t}),$$

*uniformly in $x$, dual directions $t \in \partial \mathbf{K}_\xi$ and large enough values of $r$ and $K$.*

In view of Lemma 2.3, we infer that, there exists $c_{11}$ such that, for every $\varepsilon > 0$,

$$\mathbb{P}\left(\#_{t,\delta}^{\text{mark}}(\mathfrak{t}) \geq \varepsilon \frac{|x|}{K} \,\Big|\, 0 \leftrightarrow x\right) \leq \exp(-c_{11} \varepsilon \delta |x|), \tag{2.6}$$

as soon as $K \geq K_0(\varepsilon, \delta)$, for all dual pairs $x, t$. In particular, since for *any* trunk $\mathfrak{t}$ we have a trivial a priori bound

$$N(\mathfrak{t}) \geq c_{12} \frac{|x|}{K}, \tag{2.7}$$

inequality (2.6) implies that most of the vertices of typical trunks on large enough scales are cone points. We can formulate the latter assertion as follows: There exists a constant $c_{13} > 0$ such that for every $\varepsilon > 0$ one can find $r$ and a finite scale $K_0 = K_0(\varepsilon, \delta)$ such that

$$\mathbb{P}(\#_{t,\delta}^{\text{mark}}(\mathfrak{t}) \geq \varepsilon N(\mathfrak{t}) \mid 0 \leftrightarrow x) \leq e^{-c_{13} \varepsilon |x|}, \tag{2.8}$$

for all dual pairs $x, t$ and scales $K \geq K_0$.

2.7. *Cone points of trees.* We now want to extend the previous estimate on trunks to an estimate taking also into account the branches. In order to do so, we have to slightly enlarge the opening of our forward and backward cones. Given a tree $\mathcal{T} = \mathfrak{t} \cup \mathfrak{B}$, let us say that a vertex $\mathfrak{t}_j$ of its trunk $\mathfrak{t} = (\mathfrak{t}_0, \ldots, \mathfrak{t}_j, \ldots, \mathfrak{t}_{N(\mathfrak{t})})$ is a cone point of $\mathcal{T}$ if

$$\mathcal{T} \subseteq (\mathfrak{t}_j + Y_{2\delta}^<(t)) \cup (\mathfrak{t}_j + Y_{2\delta}^>(t)) = (\mathfrak{t}_j - Y_{2\delta}^>(t)) \cup (\mathfrak{t}_j + Y_{2\delta}^>(t)). \tag{2.9}$$

Evidently, each cone point of $\mathcal{T}$ is also a cone point of $\mathfrak{t}$. The converse is, in general, not true: For $\mathfrak{t}_i \in \mathfrak{t}$, let us denote by $\mathfrak{b}(\mathfrak{t}_i)$ the branch of $\mathfrak{B}$ which is rooted at $\mathfrak{t}_i$. By Lemma 2.2, $\mathfrak{b}(\mathfrak{t}_i) = \varnothing$ for the majority of $i$'s. Some branches are, nevertheless, nonempty, so let us say that a cone point $\mathfrak{t}_j$ of $\mathfrak{t}$ is *blocked* if (2.9) does not hold. Equivalently, a cone point $\mathfrak{t}_j$ of $\mathfrak{t}$ is blocked if there exists $\mathfrak{t}_i \in \mathfrak{t}$ such that

$$\mathfrak{b}(\mathfrak{t}_i) \not\subseteq (\mathfrak{t}_j + Y_{2\delta}^{<}(t)) \cup (\mathfrak{t}_j + Y_{2\delta}^{>}(t)). \tag{2.10}$$

We claim that the number $\#_{t,\delta}^{\text{blocked}}(\mathfrak{t}, \mathfrak{B})$ of all blocked cone points of $\mathfrak{t}$ is, with overwhelming probability, small relative to the total number of all cone points of $\mathfrak{t}$. Precisely, the following uniform (in $x$, dual directions $t \in \partial \mathbf{K}_\xi$, large enough scales $K$ and admissible trunks $\mathfrak{t}$) bound holds:

Lemma 2.5. *Let an admissible trunk* $\mathfrak{t}$ *satisfy* $\#_{t,\delta}^{\text{mark}}(\mathfrak{t}) < \varepsilon N(\mathfrak{t})$. *Then there exists* $c_{14} > 0$ *such that*

$$\mathbb{P}(\#_{t,\delta}^{\text{blocked}}(\mathfrak{t}, \mathfrak{B}) \geq \varepsilon N(\mathfrak{t}) \mid \mathfrak{t}) \leq e^{-c_{14}\varepsilon |x|}. \tag{2.11}$$

Proof. As in (2.2),

$$\mathbb{P}(\mathfrak{b}(\mathfrak{t}_0) = b_0, \ldots, \mathfrak{b}(\mathfrak{t}_N) = b_N \mid \mathfrak{t}) \leq \exp\left(-(K - o_K(1)) \sum_{i=0}^{N} N(b_i)\right),$$

for each particular realization $\{b_j\}$ of $\{\mathfrak{b}(\mathfrak{t}_i)\}$. For every $i = 0, \ldots, N(\mathfrak{t})$, define $X_i$ to be the number of cone points of $\mathfrak{t}$ which are blocked by $\mathfrak{b}(\mathfrak{t}_i)$,

$$X_i = \#\{j : (2.10) \text{ holds}\}.$$

Clearly, $\#_{t,\delta}^{\text{blocked}}(\mathfrak{t}, \mathfrak{B}) \leq \sum_i X_i$. On the other hand, in view of the (cone) confinement property of cone points of $\mathcal{T}$, there exists $c_{13} > 0$ such that

$$N(\mathfrak{b}_i) \geq c_{15} X_i.$$

Therefore, exactly as in the proof of Lemma 2.2, we conclude

$$\mathbb{P}(X_1 \geq x_1, \ldots, X_N \geq x_N \mid \mathfrak{t}) \leq \exp\left(-\frac{c_{15}}{2} K \sum_{i=1}^{N} x_i\right),$$

and (2.11) follows. □

2.8. *Typical tree skeletons of* $\mathbf{C}_{0,x}$. We summarize the results of Sections 2.3–2.7 as follows: Let $x$, $t$ and $\delta$ be fixed. Given a tree skeleton $\mathcal{T} = (\mathfrak{t}, \mathfrak{B})$ of $\mathbf{C}_{0,x}$, let us denote by $\mathfrak{t}^* = (\mathfrak{t}_{i_1}, \ldots, \mathfrak{t}_{i_{N^*}})$ the collection of all $(t, \delta)$-cone points of $\mathcal{T}$. Here $N^* = N^*(\mathcal{T})$ stands for the total number of these points.

Theorem 2.1. *For every* $\delta \in (0, \frac{1}{3})$, *there exists* $r$, *a scale* $K_0$ *and two positive numbers* $\nu_3$ *and* $\nu_4$, *such that*

$$\mathbb{P}\left(N^*(\mathcal{T}) < \nu_3 \frac{|x|}{K} \,\Big|\, 0 \leftrightarrow x\right) \leq e^{-\nu_4 |x|}, \tag{2.12}$$

*uniformly in* $x$, *dual directions* $t \in \partial \mathbf{K}_\xi$ *and finite scales* $K \geq K_0$.

2.9. *Cone points of* $\mathbf{C}_{0,x}$. We say that a vertex $y \in \mathbf{C}_{0,x}$ is a $(t,\delta)$-cone point of the latter if

$$\mathbf{C}_{0,x} \subseteq \{y + Y^{<}_{3\delta}(t)\} \cup \{y + Y^{>}_{3\delta}(t)\} = \{y - Y^{>}_{3\delta}(t)\} \cup \{y + Y^{>}_{3\delta}(t)\}.$$

Notice the triple cone opening in the definition above.

If we choose $\delta$ to be too small, then it might happen that $Y^{>}_{3\delta}(t)$ does not contain axis directions, and, consequently, the above definition might be meaningless. Since, however, $\xi$ is an equivalent norm, it is not hard to see that one can choose $\delta \in (0, \frac{1}{3})$ in such a fashion that for every $t \in \partial \mathbf{K}_\xi$ there exists an axis direction $\vec{e}_i$ such that either $\vec{e}_i$ or $-\vec{e}_i$ belongs to the interior of $Y^{>}_{3\delta}(t)$. Let us fix such a value of $\delta$ for the rest of the paper.

Any cone point $\mathsf{t}_j$ of $\mathcal{T}$ is a reasonable candidate for a cone point of $\mathbf{C}_{0,x}$ itself. In fact, since by construction the Hausdorff distance between $\mathcal{T}$ and $\mathbf{C}_{0,x}$ does not exceed some $c_{16}K$, and, furthermore, away from the origin the cone $Y^{>}_{2\delta}(t)$ lies in the interior of the double cone $Y^{>}_{3\delta}(t)$, there exists a constant $c_{17}$, such that $\mathsf{t}_j$ is a cone point of $\mathbf{C}_{0,x}$ whenever

$$\mathbf{B}_{c_{17}K}(\mathsf{t}_j) \cap \mathbf{C}_{0,x} \subseteq \{\mathsf{t}_j - Y^{>}_{3\delta}(t)\} \cup \{\mathsf{t}_j + Y^{>}_{3\delta}(t)\}.$$

By the finite energy property of $\mathbb{P}$, the probability of the latter event is bounded below by a positive number $p = p(K,\delta)$ regardless of the percolation configuration on the complement $\mathbb{Z}^d \setminus \mathbf{B}_{c_{17}K}(\mathsf{t}_j)$. Therefore, for a fixed tree $\mathcal{T}$ which contains a subset $\mathsf{t}^* = (\mathsf{t}^*_1, \ldots, \mathsf{t}^*_M)$ of cone points which in addition satisfy

$$\mathbf{B}_{c_{17}K}(\mathsf{t}^*_{m+1}) \cap \mathbf{B}_{c_{17}K}(\mathsf{t}^*_m) = \varnothing,$$

the $\mathbb{P}(\cdot \mid \mathcal{T})$-conditional distribution of the number of cone points of $\mathbf{C}_{0,x}$ stochastically dominates binomial distribution $\mathsf{Bin}(M,p)$.

As before, let $\mathsf{t}^* = (t_{i_1}, \ldots, t_{i_{N^*}})$ be the collection of all cone points of $\mathcal{T}$. By Theorem 2.1, we can restrict attention to the case when $N^*(\mathcal{T}) \geq \nu_3|x|/K$. By construction, all the increments $\mathsf{t}_{i_l} - \mathsf{t}_{i_k} \in Y^{>}_{2\delta}(t)$ whenever $l > k$. Consequently,

$$\begin{aligned}\xi(\mathsf{t}_{i_l} - \mathsf{t}_{i_k}) \geq (t, \mathsf{t}_{i_l} - \mathsf{t}_{i_k})_d &= \sum_{m=k}^{l-1} (t, \mathsf{t}_{i_{m+1}} - \mathsf{t}_{i_m})_d \\ &\geq (1-\delta) \sum_{m=k}^{l-1} \xi(\mathsf{t}_{i_{m+1}} - \mathsf{t}_{i_m}) \geq (l-k)(1-\delta)K.\end{aligned}$$

Accordingly,

$$\mathbf{B}_{c_{17}K}(\mathsf{t}_{i_l}) \cap \mathbf{B}_{c_{17}K}(\mathsf{t}_{i_k}) = \varnothing \qquad \text{whenever } l - k \geq \frac{2(c_{17}+R)}{1-\delta}.$$

Let $\#^{\text{cone}}_{t,\delta}(\mathbf{C}_{0,x})$ be the number of all $(t,\delta)$-cone points of $\mathbf{C}_{0,x}$. We have proved the following:

THEOREM 2.2. *There exists $\delta \in (0, \frac{1}{3})$ and two positive numbers $\nu_5$ and $\nu_6$ such that*

$$\mathbb{P}(\#^{\text{cone}}_{t,\delta}(\mathbf{C}_{0,x}) < \nu_5|x| \mid 0 \leftrightarrow x) \le e^{-\nu_6|x|}, \tag{2.13}$$

*uniformly in $x$ and in the corresponding dual directions $t \in \partial\mathbf{K}_\xi$.*

Theorem 2.2 paves the way for an irreducible decomposition of $\mathbf{C}_{0,x}$, which we proceed to develop. Apart from (2.13), none of the estimates of Sections 2.3–2.9 is ever used in the sequel, and we shall feel free to reset the values of numerical constants $c_1$–$c_{17}$.

2.10. *Irreducible clusters.* Let $\delta$ be fixed as in Theorem 2.2, and let $t \in \partial\mathbf{K}_\xi$. A finite $R$-connected subgraph $\gamma \subset \mathbb{Z}^d$ is called backward $(t,\delta)$-irreducible, $\gamma \in \mathcal{I}^{\mathsf{b}}_t$, or just backward irreducible, if there is no ambiguity about $t$, if there exists a vertex $\mathsf{b} = \mathsf{b}(\gamma) \in \gamma$, such that

$$\gamma \subseteq \mathsf{b} + Y^{<}_{3\delta}(t) = \mathsf{b} - Y^{>}_{3\delta}(t), \tag{2.14}$$

and, in addition, $\mathsf{b}$ is the only cone point of $\gamma$.

Similarly, $\gamma$ is called forward irreducible, $\gamma \in \mathcal{I}^{\mathsf{f}}_t$, if there exists a vertex $\mathsf{f} = \mathsf{f}(\gamma) \in \gamma$ satisfying

$$\gamma \subseteq \mathsf{f} + Y^{>}_{3\delta}(t), \tag{2.15}$$

and, in addition, $\mathsf{f}$ is the only cone point of $\gamma$.

Finally, $\gamma$ is called irreducible, $\gamma \in \mathcal{I}_t$, if there exist $\mathsf{f} \neq \mathsf{b} \in \gamma$ such that both (2.14) and (2.15) are satisfied and, in addition, $\mathsf{f}$ and $\mathsf{b}$ are the only cone points of $\gamma$.

In view of Theorem 2.2, we are entitled to ignore clusters $\mathbf{C}_{0,x}$ which contain less than two $(t,\delta)$-cone points (assuming that $|x|$ is sufficiently large, of course). The remaining clusters $\mathbf{C}_{0,x}$ admit the unambiguous irreducible decomposition (1.4), where $\gamma^{\mathsf{b}}$ is backward irreducible, $\gamma^{\mathsf{f}}$ is forward irreducible, whereas the remaining clusters $\gamma_1, \ldots, \gamma_{\mathcal{N}}$ are irreducible. The number $\mathcal{N}$ of irreducible components depends on $\mathbf{C}_{0,x}$. The notation $\gamma \amalg \gamma'$ means that the two clusters are *compatible* in the sense that $\mathsf{b}(\gamma) = \mathsf{f}(\gamma')$.

**3. Probabilistic description of typical clusters.** In this section we show that typical long percolation clusters have a one-dimensional structure which complies with thermodynamical formalism of full shifts over countable alphabets.

3.1. *The representation formula.* Let $t$ and $\delta$ be fixed so that (2.13) holds. With each $(t,\delta)$-irreducible cluster $\gamma \in \mathcal{I}_t$, we associate the percolation event $\Gamma = \mathcal{W} \cap \mathcal{N} \equiv \mathcal{W}\mathcal{N}$, which is defined as follows: Let $\mathcal{E}_W(\gamma)$ be the edge set of $\gamma$. Then,

$$\mathcal{W} = \{\text{All edges of } \mathcal{E}_W(\gamma) \text{ are open}\}.$$

Accordingly, let us define the set of edges

$$\mathcal{E}_N(\gamma) = \{(u,v) : u \in \gamma\} \setminus (\mathcal{E}_W(\gamma) \cup \mathcal{E}^-(\mathsf{f}) \cup \mathcal{E}^+(\mathsf{b})), \tag{3.1}$$

where

$$\mathcal{E}^-(\mathsf{f}) \stackrel{\text{def}}{=} \{(u,\mathsf{f}) : u \in \mathsf{f} - Y^{>}_{3\delta}(t)\}$$

and

$$\mathcal{E}^+(\mathsf{b}) = \{(\mathsf{b},v) : v \in \mathsf{b} + Y^{>}_{3\delta}(t)\},$$

and, of course, $\mathsf{f} = \mathsf{f}(\gamma)$ and $\mathsf{b} = \mathsf{b}(\gamma)$. Then,

$$\mathcal{N} = \{\text{All edges of } \mathcal{E}_N(\gamma) \text{ are closed}\}.$$

Notice that we refrain from closing edges from $\mathcal{E}^-(\mathsf{f})$ and $\mathcal{E}^+(\mathsf{b})$ in (3.1) in order to be compatible with the decomposition (1.4) of $\mathbf{C}_{0,x}$.

The definition of percolation events $\Gamma^{\mathsf{b}} = \mathcal{W}^{\mathsf{b}}\mathcal{N}^{\mathsf{b}}$ and $\Gamma^{\mathsf{f}} = \mathcal{W}^{\mathsf{f}}\mathcal{N}^{\mathsf{f}}$, which correspond to backward and, respectively, forward irreducible clusters, is completely similar.

In this notation (1.4) gives rise [up to factors of order $1 + o(e^{-\nu_6|x|})$] to the following representation formula:

$$\mathbb{P}(0 \leftrightarrow x) = \sum_{\gamma^{\mathsf{b}} \ni 0} \sum_{\gamma^{\mathsf{f}} \ni x} \sum_{n} \sum^{*}_{\gamma_1,\dots,\gamma_n} \mathbb{P}(\Gamma^{\mathsf{b}}\Gamma_1 \cdots \Gamma_n \Gamma^{\mathsf{f}}), \tag{3.2}$$

where the last sum is over all strings of irreducible compatible clusters $\underline{\gamma} = (\gamma_1,\dots,\gamma_n)$ (and with respect to associated percolation events $\Gamma_1,\dots,\Gamma_n$).

It is convenient to consider all the clusters in question modulo their $\mathbb{Z}^d$-shifts. In this way any collection $\gamma^{\mathsf{b}}, \underline{\gamma} = (\gamma_1,\dots,\gamma_n), \gamma^{\mathsf{f}}$ of, respectively, backward, string of irreducible, and forward clusters can be patched together. Introducing an additional notation

$$V(\gamma) \stackrel{\text{def}}{=} \mathsf{b}(\gamma) - \mathsf{f}(\gamma),$$

for the displacement along an irreducible cluster $\gamma$, and the half-spaces

$$\mathcal{H}_0^{t,+} \stackrel{\text{def}}{=} \{y \in \mathbb{Z}^d : (t,y)_d \geq 0\}$$

and

$$\mathcal{H}_x^{t,-} \stackrel{\text{def}}{=} \{y \in \mathbb{Z}^d : (t,x)_d \geq (t,y)_d\},$$

we rewrite (3.2) as

$$
\begin{aligned}
&e^{\xi(x)}\mathbb{P}(0\leftrightarrow x)\\
&\quad=\sum_{y\in\mathcal{H}_0^{t,+}}\sum_{z\in\mathcal{H}_x^{t,-}}\sum_{\gamma^{\mathsf{b}}\ni -y}\sum_{\gamma^{\mathsf{f}}\ni x-z}\mathbb{P}(\Gamma^{\mathsf{b}})\mathbb{P}(\Gamma^{\mathsf{f}})e^{(t,x-z+y)_d}\\
&\qquad\times\sum_n\sum_{V(\gamma_1)+\cdots+V(\gamma_n)=z-y}\frac{\mathbb{P}(\Gamma^{\mathsf{b}}\Gamma_1\cdots\Gamma_n\Gamma^{\mathsf{f}})}{\mathbb{P}(\Gamma^{\mathsf{b}})\mathbb{P}(\Gamma^{\mathsf{f}})}e^{(t,z-y)_d},
\end{aligned}
\tag{3.3}
$$

where the shifts of backward and forward clusters are normalized in such a way that $\mathsf{b}(\gamma^{\mathsf{b}})=\mathsf{f}(\gamma^{\mathsf{f}})=0$.

By Theorem 2.2, the above equality holds up to factors of order $1+o(e^{-\nu_6|x|})$.

3.2. *Decomposition of probabilities.* We shall decompose the probabilities in (3.3) as follows: First of all set

$$
\Gamma^{\mathsf{b}}\Gamma_1\cdots\Gamma_n\Gamma^{\mathsf{f}}=\mathcal{W}^{\mathsf{b}}\mathcal{N}^{\mathsf{b}}\mathcal{W}_1\mathcal{N}_1\cdots\mathcal{W}_n\mathcal{N}_n\mathcal{W}^{\mathsf{f}}\mathcal{N}^{\mathsf{f}}\stackrel{\text{def}}{=}\mathcal{W}^{\mathsf{b}}\mathcal{N}^{\mathsf{b}}\mathcal{W}\mathcal{N},
$$

where $\mathcal{W}=\mathcal{W}_1\cdots\mathcal{W}^{\mathsf{f}}$ and $\mathcal{N}=\mathcal{N}_1\cdots\mathcal{N}^{\mathsf{f}}$.

At this point we shall rely on a specific decoupling property of FK measures[7]: namely,

$$
\mathbb{P}(\mathcal{W}^{\mathsf{b}}\mathcal{W}\mid\mathcal{N}^{\mathsf{b}}\mathcal{N})=\mathbb{P}^{\mathsf{f}}_{\gamma^{\mathsf{b}}\amalg\gamma_1\amalg\cdots\amalg\gamma_n\amalg,\gamma^{\mathsf{f}}}(\mathcal{W}^{\mathsf{b}}\mathcal{W})=\mathbb{P}^{\mathsf{f}}_{\gamma^{\mathsf{b}}}(\mathcal{W}^{\mathsf{b}})\mathbb{P}^{\mathsf{f}}_{\gamma^{\mathsf{f}}}(\mathcal{W}^{\mathsf{f}})\prod_1^n\mathbb{P}^{\mathsf{f}}_{\gamma_k}(\mathcal{W}_k).
$$

In the same fashion,

$$
\mathbb{P}(\Gamma^{\mathsf{b}})\mathbb{P}(\Gamma^{\mathsf{f}})=\mathbb{P}^{\mathsf{f}}_{\gamma^{\mathsf{b}}}(\mathcal{W}^{\mathsf{b}})\mathbb{P}^{\mathsf{f}}_{\gamma^{\mathsf{f}}}(\mathcal{W}^{\mathsf{f}})\mathbb{P}(\mathcal{N}^{\mathsf{b}})\mathbb{P}(\mathcal{N}^{\mathsf{f}}).
$$

Therefore,

$$
\frac{\mathbb{P}(\Gamma^{\mathsf{b}}\Gamma_1\cdots\Gamma_n\Gamma^{\mathsf{f}})}{\mathbb{P}(\Gamma^{\mathsf{b}})\mathbb{P}(\Gamma^{\mathsf{f}})}=\prod_{k=1}^n\mathbb{P}^{\mathsf{f}}_{\gamma_k}(\mathcal{W}_k)\frac{\mathbb{P}(\mathcal{N}^{\mathsf{b}}\mathcal{N}_1\cdots\mathcal{N}_n\mathcal{N}^{\mathsf{f}})}{\mathbb{P}(\mathcal{N}^{\mathsf{b}})\mathbb{P}(\mathcal{N}^{\mathsf{f}})}.
$$

Finally, we decompose the above fraction as

$$
\begin{aligned}
&\frac{\mathbb{P}(\mathcal{N}^{\mathsf{b}}\mathcal{N}_1\cdots\mathcal{N}_n\mathcal{N}^{\mathsf{f}})}{\mathbb{P}(\mathcal{N}^{\mathsf{b}})\mathbb{P}(\mathcal{N}^{\mathsf{f}})}\\
&\quad=\frac{\mathbb{P}(\mathcal{N}^{\mathsf{b}}\mid\mathcal{N}_1\cdots\mathcal{N}_n\mathcal{N}^{\mathsf{f}})}{\mathbb{P}(\mathcal{N}^{\mathsf{b}})}\,\mathbb{P}(\mathcal{N}_1\mid\mathcal{N}_2\cdots\mathcal{N}_n\mathcal{N}^{\mathsf{f}})\cdots\mathbb{P}(\mathcal{N}_n\mid\mathcal{N}^{\mathsf{f}}).
\end{aligned}
\tag{3.4}
$$

[7]Notice that the rest of our arguments depends only on positive association, finite energy and assumption (1.2).

Alternatively, we could have decomposed them as

$$\begin{aligned}(3.5)\qquad &\frac{\mathbb{P}(\mathcal{N}^{\mathsf{b}}\mathcal{N}_1\cdots\mathcal{N}_n\mathcal{N}^{\mathsf{f}})}{\mathbb{P}(\mathcal{N}^{\mathsf{b}})\mathbb{P}(\mathcal{N}^{\mathsf{f}})}\\ &\quad=\frac{\mathbb{P}(\mathcal{N}^{\mathsf{f}}\,|\,\mathcal{N}_n\cdots\mathcal{N}_1\mathcal{N}^{\mathsf{b}})}{\mathbb{P}(\mathcal{N}^{\mathsf{f}})}\,\mathbb{P}(\mathcal{N}_n\,|\,\mathcal{N}_{n-1}\cdots\mathcal{N}_1\mathcal{N}^{\mathsf{b}})\cdots\mathbb{P}(\mathcal{N}_1\,|\,\mathcal{N}^{\mathsf{b}}).\end{aligned}$$

A comparison between (3.4) and (3.5) enables to formulate results in a particularly convenient symmetric form.

3.3. *Reduction to one-dimensional thermodynamics.* With clusters $\gamma^{\mathsf{b}}$ and $\gamma^{\mathsf{f}}$ fixed, define the function

$$(3.6)\qquad g(\gamma_1,\dots,\gamma_n)=g(\gamma_1,\dots,\gamma_n\,|\,\gamma^{\mathsf{b}},\gamma^{\mathsf{f}})\stackrel{\mathrm{def}}{=}\frac{\mathbb{P}(\mathcal{N}^{\mathsf{b}}\,|\,\mathcal{N}_1\cdots\mathcal{N}_n\mathcal{N}^{\mathsf{f}})}{\mathbb{P}(\mathcal{N}^{\mathsf{b}})},$$

and a potential

$$(3.7)\quad e^{\psi_t(\gamma_1,\dots,\gamma_n)}=e^{\psi_t(\gamma_1,\dots,\gamma_n|\gamma^{\mathsf{f}})}\stackrel{\mathrm{def}}{=}e^{(t,V(\gamma_1))_d}\mathbb{P}^{\mathrm{f}}_{\gamma_1}(\mathcal{W}_1)\mathbb{P}(\mathcal{N}_1\,|\,\mathcal{N}_2\cdots\mathcal{N}_n\mathcal{N}^{\mathsf{f}}).$$

Both the function $g(\cdot)$ and the potential $\psi_t(\cdot)$ act on any string $\underline{\gamma}$ of length $n$ of $(\mathcal{I}_t)$ irreducible clusters for $n=1,2,\dots$. For any $n\in\mathbb{N}$, define a measure $\mu^t_n(\cdot\,|\,\gamma^{\mathsf{b}},\gamma^{\mathsf{f}})$ on the set $\Xi_n$ of strings of $n$-irreducible clusters

$$\mu^t_n(\gamma_1,\dots,\gamma_n\,|\,\gamma^{\mathsf{b}},\gamma^{\mathsf{f}})\stackrel{\mathrm{def}}{=}e^{\psi_t(\gamma_1,\gamma_2,\dots,\gamma_n)+\psi_t(\gamma_2,\dots,\gamma_n)+\cdots+\psi_t(\gamma_n)}g(\gamma_1,\dots,\gamma_n).$$

In this notation the last sum in (3.3) is just

$$(3.8)\qquad \mu^t_n(V(\gamma_1)+\cdots+V(\gamma_n)=z-y\,|\,\gamma^{\mathsf{b}},\gamma^{\mathsf{f}}).$$

It is easy to see that Theorem 2.2 implies that there exists $\nu_7>0$, such that

$$(3.9)\qquad \sum_{\gamma^{\mathsf{b}}\ni -u}\mathbb{P}(\Gamma^{\mathsf{b}})e^{(t,u)_d}\le e^{-\nu_7|u|}\quad\text{and}\quad\sum_{\gamma^{\mathsf{f}}\ni u}\mathbb{P}(\Gamma^{\mathsf{f}})e^{(t,u)_d}\le e^{-\nu_7|u|},$$

uniformly in $u\in\mathbb{Z}^d$. Consequently, the main contribution to the right-hand side of (3.3) should come from the terms (3.8), with $z-y$ being close to $x$ and with the number of steps $n$ being close to optimal. The local limit analysis of measures $\mu^t_n(\cdot\,|\,\gamma^{\mathsf{b}},\gamma^{\mathsf{f}})$ boils down to a study of the analytic spectral properties of the Ruelle type operator

$$(3.10)\qquad \mathfrak{L}_t f(\underline{\gamma})\stackrel{\mathrm{def}}{=}\sum_{\gamma\ \text{irreducible}}e^{\psi_t(\gamma,\underline{\gamma})}f(\gamma,\underline{\gamma}).$$

In this respect we shall be able to rely on a general theory developed in [9, 10] as soon as we check appropriate regularity properties of $g(\cdot)$ and $\psi_t(\cdot)$. This is the content of the next section.

3.4. *Regularity properties of* $g(\cdot)$ *and* $\psi_t(\cdot)$. For two strings of irreducible clusters $\underline{\gamma} = (\gamma_1, \ldots, \gamma_n)$ and $\underline{\lambda} = (\lambda_1, \ldots, \lambda_m)$, define

$$\mathfrak{i}(\underline{\gamma}, \underline{\lambda}) \stackrel{\text{def}}{=} \min\{i : \lambda_i \neq \gamma_i\} \wedge n \wedge m.$$

THEOREM 3.1. *There exist* $\theta < 1$ *and two positive constants* $\nu_8$ *and* $\nu_9$, *such that the following inequalities hold uniformly in* $t \in \partial\mathbf{K}_\xi$, $t$*-irreducible clusters* $\gamma$, *strings of* $t$*-irreducible clusters* $\underline{\gamma}, \underline{\lambda}$, $t$*-backward irreducible clusters* $\gamma^{\mathsf{b}}$ *and pairs of* $t$*-forward irreducible clusters* $\gamma^{\mathsf{f}}$, $\tilde{\gamma}^{\mathsf{f}}$:

$$\nu_8 \leq g(\underline{\gamma} \mid \gamma^{\mathsf{b}}, \gamma^{\mathsf{f}}) \leq \frac{1}{\nu_8}. \tag{3.11}$$

*Furthermore,*

$$|g(\gamma, \underline{\gamma} \mid \gamma^{\mathsf{b}}, \gamma^{\mathsf{f}}) - g(\gamma, \underline{\lambda} \mid \gamma^{\mathsf{b}}, \tilde{\gamma}^{\mathsf{f}})| \leq \nu_9 \theta^{\mathfrak{i}(\underline{\gamma}, \underline{\lambda})}, \tag{3.12}$$

*as well as*

$$|\psi_t(\gamma, \underline{\gamma} \mid \gamma^{\mathsf{f}}) - \psi_t(\gamma, \underline{\lambda} \mid \tilde{\gamma}^{\mathsf{f}})| \leq \nu_9 \theta^{\mathfrak{i}(\underline{\gamma}, \underline{\lambda})}. \tag{3.13}$$

PROOF. Given a set of bonds $A$, define

$$\mathcal{N}(A) = \bigcap_{e \in A} \{\omega(e) = 0\}.$$

Thanks to the cone confinement property of irreducible clusters, all three claims (3.11), (3.12) and (3.13) of Theorem 3.1 are straightforward consequences of the following

PROPOSITION 3.1. *Assume that* $\beta < \beta_c^2$. *Then, there exist a constant* $c_1$, *such that*

$$\frac{\mathbb{P}(\mathcal{N}(A) \mid \mathcal{N}(B))}{\mathbb{P}(\mathcal{N}(A) \mid \mathcal{N}(C))} \leq \exp\left\{c_1 \sum_{e \in A} e^{-\nu_1 \mathrm{d}(e, B\Delta C)}\right\}, \tag{3.14}$$

*uniformly in (finite) sets of bond* $A$, $B$ *and* $C$.

It is enough to consider $C = B \vee D$. Then, for every bond $e$, the FKG property of $\mathbb{P}$ implies

$$0 \leq \mathbb{P}(\omega(e) = 0 \mid \mathcal{N}(B \vee D)) - \mathbb{P}(\omega(e) = 0 \mid \mathcal{N}(B)).$$

However, a repeated application of the FKG property of $\mathbb{P}$ implies

$$\begin{aligned}\mathbb{P}(\omega(e) = 0 \mid \mathcal{N}(B)) &\geq \mathbb{P}(\omega(e) = 0; e \not\leftrightarrow D \mid \mathcal{N}(B)) \\ &\geq \mathbb{P}(\omega(e) = 0 \mid e \not\leftrightarrow D; \mathcal{N}(B))\mathbb{P}(e \not\leftrightarrow D) \\ &\geq \mathbb{P}(\omega(e) = 0 \mid \mathcal{N}(B \vee D))(1 - \nu_0 e^{-\nu_1 \mathrm{d}(e, D)}),\end{aligned}$$

where we have also used exponential decay of connectivities on the last step. As a result,

$$\begin{aligned}0 &\le \mathbb{P}(\omega(e)=0 \,|\, \mathcal{N}(B\vee D)) - \mathbb{P}(\omega(e)=0\,|\,\mathcal{N}(B)) \\ &\le \nu_0 e^{-\nu_1 \mathrm{d}(e,D)}\mathbb{P}(\omega(e)=0\,|\,\mathcal{N}(B\vee D)).\end{aligned}$$

Thus,

$$1 \ge \frac{\mathbb{P}(\omega(e)=0\,|\,\mathcal{N}(B))}{\mathbb{P}(\omega(e)=0\,|\,\mathcal{N}(B\vee D))} \ge 1-\nu_0 e^{-\nu_1 \mathrm{d}(e,D)},$$

and (3.14) follows for any singleton $A=\{e\}$.

The general case reduces to the singleton one. Let us number the bonds of $A$ as $A=\{e_1,\dots,e_r\}$. After that set $A_k=\{e_1,\dots,e_k\}$. Clearly,

$$\frac{\mathbb{P}(\mathcal{N}(A)\,|\,\mathcal{N}(B))}{\mathbb{P}(\mathcal{N}(A)\,|\,\mathcal{N}(B\vee D))} = \prod_0^{r-1}\frac{\mathbb{P}(\omega(e_{k+1})=0\,|\,\mathcal{N}(A_k\vee B))}{\mathbb{P}(\omega(e_{k+1})=0\,|\,\mathcal{N}(A_k\vee B\vee D))},$$

and (3.14) follows. □

## 4. Implications.

We are ready now to explain all our main results.

4.1. *Ornstein–Zernike formula for two point functions.* Theorem 3.1 implies that the local limit analysis of (3.8) falls into the framework of the general local limit theory developed in Sections 4 and 5 of [9]. In particular, only terms with

$$|(z-y)-x| \ll |x| \tag{4.1}$$

really contribute to the right-hand side of (3.3). Furthermore, there exist a function $F$ on the set of all backward and forward irreducible clusters, which is bounded above and below,

$$\frac{1}{c_2} \le F(\cdot) \le c_2,$$

for some $c_2 \ge 1$, and a locally analytic positive function $\hat{\Psi}$ defined on a neighborhood of $\vec{n}_x$ in $\mathbb{S}^{d-1}$, so that

$$\begin{aligned}&\sum_n \mu_n^t(V(\gamma_1)+\cdots+V(\gamma_n)=z-y\,|\,\gamma^{\mathsf{b}},\gamma^{\mathsf{f}}) \\ &\quad = \frac{\hat{\Psi}(\vec{n}_x)}{|x|^{(d-1)/2}} e^{-\xi(z-y)+(t,z-y)_d} F(\gamma^{\mathsf{b}})F(\gamma^{\mathsf{f}})(1+o(1)).\end{aligned} \tag{4.2}$$

Since $-\xi(z-y)+(t,z-y)_d \le 0$, (3.9) implies that the range of $y$ and $z$ which essentially contributes to the right-hand side of (3.3) can be further shrunk to $|y|, |z-x| \le c_3 \log|x|$ for some $c_3$ large enough. For such $y$ and $z$, however,

$$\xi(z-y) - (t, z-y)_d = \xi(z-y) - \xi(x) - (t, (z-y)-x)_d$$
$$= O\biggl(\frac{|(z-y)-x|^2}{|x|}\biggr) = o(1)$$

uniformly in $|x| \to \infty$. As a result, the assertion of Theorem A follows with

$$\Psi(\vec{n}_x) = \biggl(\sum_{u \in \mathcal{H}_0^{t,+}} \sum_{\gamma^{\mathsf{f}} \ni u} \mathbb{P}(\Gamma^{\mathsf{f}}) e^{(t,u)_d} F(\gamma^{\mathsf{f}})\biggr)^2 \cdot \hat{\Psi}(\vec{n}_x),$$

where we have relied on lattice symmetries of the model, and in particular, on the central symmetry employed in the definition of backward and forward irreducible clusters, on the last step.

4.2. *Exit from boxes*: *Proof of Corollary* 1.1. The proof of Corollary 1.1 is very similar to the derivation of Ornstein–Zernike asymptotics, so we only sketch the argument.

First of all, observe that $\mathbf{K}_\xi$ inherits all the symmetries of $\mathbf{J}$. Therefore, its strict convexity implies that the minimum of $\xi(x)$ on $\{x \in \mathbb{R}^d : |x|_\infty = 1\}$ is attained exactly when $x = \pm\vec{e}_i$, $i = 1, \dots, d$.

Next, one has the following obvious lower bound on the probability of exiting the box $\Lambda_N$: let $x^N = (N+1)\vec{e}_1$, then

$$\mathbb{P}(0 \leftrightarrow \mathbb{Z}^d \setminus \Lambda_N) \ge \mathbb{P}(0 \leftrightarrow x^N) \ge CN^{-(d-1)/2} e^{-N\xi(\vec{e}_1)},$$

for some $C = C(\beta, \mathbf{J}, d) > 0$. Let $t \in \partial\mathbf{K}_\xi$ be dual to $x^N$ (notice that $t$ is a multiple of $\vec{e}_1$). Using the previous lower bound, one obtains that

$$\mathbb{P}(0 \leftrightarrow \mathbb{Z}^d \setminus \Lambda_N) = 2d\mathbb{P}(0 \leftrightarrow Y_{2\delta}^{>}(t) \setminus \Lambda_N)(1 + o(1)),$$

so the problem is reduced to the analysis of the probability in the right-hand side.

As before, up to factors of order $1 + o(e^{-\nu_6 N})$, clusters connecting 0 to the exterior of the box $\Lambda_N$, in the general direction specified by the forward cone $Y_{3\delta}^{>}(t)$, admit a decomposition of the type (3.3),

$$e^{\xi(x^N)} \mathbb{P}(0 \leftrightarrow Y_{3\delta}^{>}(t) \setminus \Lambda_N)$$
$$= \sum_{y,z \in \Lambda_N} \sum_{\gamma^{\mathsf{b}} \ni -y} \sum_{\gamma^{\mathsf{f}} : \mathcal{R}(\gamma) \ge (\vec{e}_1, x^N - z)_d} \mathbb{P}(\Gamma^{\mathsf{b}}) \mathbb{P}(\Gamma^{\mathsf{f}}) e^{(t, x^N - z + y)_d}$$
$$\times \sum_n \sum_{V(\gamma_1) + \cdots + V(\gamma_n) = z - y} \frac{\mathbb{P}(\Gamma^{\mathsf{b}} \Gamma_1 \cdots \Gamma_n \Gamma^{\mathsf{f}})}{\mathbb{P}(\Gamma^{\mathsf{b}}) \mathbb{P}(\Gamma^{\mathsf{f}})} e^{(t, z-y)_d},$$

where the horizontal range of a cluster $C$ is defined as

$$\mathcal{R}(C) \stackrel{\text{def}}{=} \sup_{u\in C}(\vec{e}_1, u)_d - \inf_{u\in C}(\vec{e}_1, u)_d.$$

It is convenient to write $z = (x_1^N - \hat{z}_1, \hat{z}_\perp)$, with $\hat{z}_1 \in \mathbb{Z}$ and $\hat{z}_\perp \in \mathbb{Z}^{d-1}$. Reasoning as in Section 4.1, we see that one can restrict attention to points $y$ and $z$ such that $|y| \leq c_3 \log N$ and $\hat{z}_1 \leq c_3 \log N$, for some sufficiently large constant $c_3$. In this case, writing $\hat{z} = (N+1, \hat{z}_\perp)$, we obtain that

$$\begin{aligned}\sum_n \mu_n^t(V(\gamma_1) + \cdots + V(\gamma_n) = z - y \mid \gamma^{\mathsf{b}}, \gamma^{\mathsf{f}})\\ = \frac{\hat{\Psi}(\vec{n}_{\hat{z}})}{|\hat{z}|^{(d-1)/2}} e^{-\xi(\hat{z})+\xi(\vec{e}_1)N} F(\gamma^{\mathsf{b}})F(\gamma^{\mathsf{f}})(1+o(1)).\end{aligned}$$

Rough estimates show now that the contribution of terms with $|\hat{z}_\perp| > N^{1/2+\eta}$, $\eta > 0$, is negligible. For the remaining terms, we have

$$\begin{aligned}\sum_{|\hat{z}_\perp| \leq N^{1/2+\eta}} \frac{\hat{\Psi}(\vec{n}_{\hat{z}})}{|\hat{z}|^{(d-1)/2}} e^{-\xi(\hat{z})+\xi(\vec{e}_1)N}\\ = \frac{\hat{\Psi}(\vec{e}_1)}{N^{(d-1)/2}} \sum_{|\hat{z}_\perp| \leq N^{1/2+\eta}} e^{-N(\hat{z}_\perp, \mathrm{d}^2\xi|_{\vec{e}_1}\hat{z}_\perp)_d/2}(1+o(1))\\ = \hat{\Psi}(\vec{e}_1) \int d^{d-1}u \; e^{-(u, \mathrm{d}^2\xi|_{\vec{e}_1}u)_d/2}(1+o(1)),\end{aligned}$$

which completes the proof.

4.3. *Geometry of Wulff shapes and equi-decay profiles.* By the relation between $\xi$ and $\mathbf{K}_\xi$,

$$\xi(x) = \max_{t\in K}(t, x)_d,$$

the Wulff shape $\mathbf{K}_\xi$ can be alternatively defined as the closure of the convergence set

$$\left\{ t : \sum_x \mathbb{P}(0 \leftrightarrow x) e^{(t,x)_d} < \infty \right\}.$$

Consequently (see Section 3.3 in [9]), given any $t \in \partial\mathbf{K}_\xi$, there exists $\varepsilon > 0$ such that the portion of the boundary $\partial\mathbf{K}_\xi$ inside the $\varepsilon$-ball $t + B_\varepsilon$ around $t$ can be described as the level set

$$\partial\mathbf{K}_\xi \cap (t + B_\varepsilon) = \{t + s : \rho(\mathfrak{L}_{t,s}) = 1\},$$

where $\rho(\mathfrak{L}_{t,s})$ is the spectral radius of the operator

$$\mathfrak{L}_{t,s} f(\underline{\gamma}) = \sum_{\gamma \in \mathcal{I}_t} e^{\psi_t(\gamma, \underline{\gamma}) + (s, V(\gamma))_d} f(\gamma, \underline{\gamma}).$$

Of course, $\mathfrak{L}_{t,s}$ is a perturbation of the operator $\mathfrak{L}_t$ defined in (3.10). In view of Theorem 3.1, the operator $\mathfrak{L}_t$ satisfies all the assumptions of general Theorem 4.1 in [9]. In particular, $\rho(\mathfrak{L}_t)$ is a nondegenerate eigenvalue of $\mathfrak{L}_t$ and the Fredholm spectrum of $\mathfrak{L}_t$ lies strictly inside the $\mathbb{C}$-ball of radius $\rho(\mathfrak{L}_t)$. Local analyticity of $\partial \mathbf{K}_\xi$ around $t$ follows then from analytic perturbation theory of nondegenerate eigenvalues. Positivity of Gaussian curvature follows from the conditional variance argument. We refer to Sections 5.4 and 5.5 of [9] for the corresponding details.

4.4. *Invariance principle for long connected clusters.* In view of the decomposition (1.4) and the regularity properties of the potentials as stated in Theorem 3.1, the proof of Section 3 of [14] just goes through.

4.5. *Nearest neighbor Potts models in the phase transition regime.* Using the standard coupling of the Potts model with $\vec{n}$ Dobrushin boundary condition and the associated random-cluster measure [21], it is easy to check that in any pair of configurations $(\omega, \sigma)$ contributing to the joint measure one has the following property: All bonds dual to those belonging to the interface of $\sigma$ are closed in $\omega$ (indeed, they join sites with spins taking different values). By duality, this means that all the bonds belonging to the interface are open in $\omega^*$. Therefore, the distribution of the bonds in the interface is actually stochastically dominated by the distribution of the bonds in the cluster $C_{x_L^N, x_R^N}$ connecting the dual sites $x_L^N \stackrel{\text{def}}{=} (-N - \frac{1}{2}, -\frac{1}{2})$ and $x_R^N \stackrel{\text{def}}{=} (N + \frac{1}{2}, \frac{1}{2})$ under the measure

$$\mathbb{P}_\Lambda^{\mathsf{f}}(\cdot \mid x_L^N \leftrightarrow x_R^N)$$

at inverse temperature $\beta^*$ (i.e., in the subcritical regime). The theorem would therefore follow if we can extend Theorem C to this finite-volume situation. As in [14], the only required change to the proof is to redefine the two extremal irreducible pieces $\gamma^{\mathsf{b}}$ and $\gamma^{\mathsf{f}}$: one introduces two new extremal pieces $\bar{\gamma}^{\mathsf{b}} \stackrel{\text{def}}{=} \gamma^{\mathsf{b}} \amalg \gamma_1 \amalg \cdots \amalg \gamma_{\lfloor (\log N)^2 \rfloor}$ and $\bar{\gamma}^{\mathsf{f}} \stackrel{\text{def}}{=} \gamma_{\mathcal{N} - \lfloor (\log N)^2 \rfloor} \amalg \cdots \amalg \gamma_{\mathcal{N}} \amalg \gamma^{\mathsf{f}}$. The size of $\bar{\gamma}^{\mathsf{b}}$ and $\bar{\gamma}^{\mathsf{f}}$ can be assumed to be at most $O(\log N)^4$ and therefore does not affect the limiting behavior, but is still large enough to keep the boundary effects under control. We refer to the detailed discussion in [14] (dealing with the case of the Ising model).

## APPENDIX: BASIC ASSUMPTION IN TWO DIMENSIONS

In two dimensions our assumption (1.2) holds once the inverse correlation length is positive.

PROPOSITION A.1. *Let $\{\mathbf{J}, \beta, q\}$ be fixed and assume that there is a unique $\{\mathbf{J}, \beta, q\}$-infinite volume measure $\mathbb{P}$ which, in addition, has exponential decay of connectivities,*

$$\mathbb{P}(0 \leftrightarrow x) \le e^{-c_1 |x|}.$$

*Then assumption (1.2) is satisfied.*

The above proposition has been established in a somewhat more general context in [3]. Here we sketch a simplified proof which was kindly explained to us by Reda Messikh (see also [11]) with a reference to the original idea by David Barbato. In order to facilitate the notation, we shall consider the nearest-neighbor case, so that the assumption on positivity of inverse correlation length is equivalent to the assumption on positivity of surface tension in the dual model. An adjustment to the general finite-range case is straightforward.

STEP 1. Define

$$\eta(x) = I_{\{x \leftrightarrow \partial \Lambda_N\}}.$$

Then, $\sum_{x \in \Lambda_N} \eta(x)$ is just the size of the boundary cluster. We claim that, for every $\delta > 0$,

$$\limsup_{N \to \infty} \frac{1}{N} \mathbb{P}^{\mathrm{w}}_{\Lambda_N} \left( \frac{1}{N^2} \sum_{x \in \Lambda_N} \eta(x) \ge \delta \right) = -\infty. \tag{A.1}$$

The super-exponential bound (A.1) enables to restrict attention to arbitrary small boundary cluster densities. In fact, a stronger claim is true:

$$\limsup_{N \to \infty} \frac{1}{N^2} \mathbb{P}^{\mathrm{w}}_{\Lambda_N} \left( \frac{1}{|\Lambda_N|} \sum_{x \in \Lambda_N} \eta(x) \ge \delta \right) < 0, \tag{A.2}$$

for every $\delta > 0$ fixed. (A.2) follows from the Hoeffding inequality: First of all,

$$\lim_{M \to \infty} \mathbb{P}^{\mathrm{w}}_{\Lambda_M} (0 \leftrightarrow \partial \Lambda_M) = 0,$$

by the unicity of the infinite volume measure $\mathbb{P}$. Thus, for every $\delta > 0$, one can choose $M = M(\delta) < \infty$, such that

$$\mathbb{E}^{\mathrm{w}}_{\Lambda_M} \left( \frac{1}{|\Lambda_M|} \sum_{x \in \Lambda_M} I_{\{x \leftrightarrow \Lambda_M\}} \right) \le \frac{\delta}{2}. \tag{A.3}$$

Without loss of generality, we may assume that $M(\delta)$ divides $N$, $N = (2K + 1)M$. Thus,

$$\Lambda_N = \bigvee_{x \in \Lambda_K} (2Mx + \Lambda_M).$$

Since

$$\frac{1}{|\Lambda_N|} \sum_{x \in \Lambda_N} \eta(x) \le \frac{1}{|\Lambda_K|} \sum_{x \in \Lambda_K} \left( \frac{1}{|\Lambda_M|} \sum_{y \in 2Mx + \Lambda_M} I_{\{y \leftrightarrow \partial(2Mx + \Lambda_M)\}} \right),$$

it remains to use (A.3) and positive association of random cluster measures.

Step 2. Next one constructs a minimal section type splitting of the event $\{0 \leftrightarrow \partial\Lambda_{N/2}\}$. In the sequel we shall use the notation

$$\Lambda_N^l \stackrel{\text{def}}{=} \Lambda_N \setminus \Lambda_{N-l-1}.$$

Fix a small $\delta > 0$. By the preceding step, we may assume that

$$\sum_{x \in \Lambda_N^{N/2}} \eta(x) < \frac{\delta N^2}{2}. \tag{A.4}$$

In particular, at least one square layer (below $|\cdot|$ is the maximum norm) $L_N^l \stackrel{\text{def}}{=} \{x : |x| = N - l\}$, $l = 1, \ldots, N/2$, has less that $\delta N$ sites connected to the outer boundary $\partial\Lambda_N$. It is convenient to decouple the corresponding events from the event $\{0 \leftrightarrow \partial\Lambda_{N/2}\}$ as follows: Consider random variables

$$\omega(x) = \begin{cases} 1, & \text{if } x \stackrel{\Lambda_N^{|x|}}{\longleftrightarrow} \partial\Lambda_N, \\ 0, & \text{otherwise,} \end{cases}$$

and set

$$\Omega(l) = \sum_{|x| = N-l} \omega(x).$$

Of course, $\omega(x) \le \eta(x)$, and

$$l^* \stackrel{\text{def}}{=} \min\{l : \Omega(l) < \delta N\} \in \left\{1, \ldots, \frac{N}{2} - 1\right\},$$

for any bond configuration satisfying (A.4). Define the event $A_l = \{l^* = l\}$. Then, ignoring configurations which violate (A.4), we arrive to the following decomposition:

$$\{0 \leftrightarrow \partial\Lambda_{N/2}\} = \bigvee_{l=1}^{N/2-1} \{0 \leftrightarrow \partial\Lambda_{N/2}; A_l\}.$$

Finally, note that $A_l$ depends only on the percolation configuration on bonds $(x,y)$ with both end-points belonging to $\Lambda_N^l$. For each $l = 1,\dots,N/2-1$, the latter set of bonds is disjoint from the set of bonds with both end-points in $\Lambda_{N/2}$.

STEP 3. Assume $A_l$ and fix a bond configuration $\phi_l$ on $\Lambda_N^l$ and a bond configuration $\psi_l$ on $\Lambda_{N-l-1}$. Since there are less than $\delta N$ sites of the layer $L_{N-l}$ which are connected (in $\phi_l$) to $\partial\Lambda_N$, the finite energy condition implies that

$$\mathbb{P}_{\Lambda_N}^{\mathrm{w}}(\phi_l;\psi_l) \le e^{c_2\delta N}\mathbb{P}_{\Lambda_N}^{\mathrm{w}}(\phi_l;\psi_l;\partial\Lambda_{N/2} \not\leftrightarrow \partial\Lambda_N)$$

for some universal constant $c_2$. It follows that, for each $l$,

$$\mathbb{P}_{\Lambda_N}^{\mathrm{w}}(0 \leftrightarrow \partial\Lambda_{N/2}; A_l) \le e^{c_2\delta N}\mathbb{P}_{\Lambda_N}^{\mathrm{w}}(0 \leftrightarrow \partial\Lambda_{N/2}; \Lambda_{N/2} \not\leftrightarrow \partial\Lambda_N; A_l).$$

However,

$$\mathbb{P}_{\Lambda_N}^{\mathrm{w}}(0 \leftrightarrow \partial\Lambda_{N/2}; \Lambda_{N/2} \not\leftrightarrow \partial\Lambda_N) \le \mathbb{P}(0 \leftrightarrow \partial\Lambda_{N/2})$$

by a standard FKG decoupling argument.

**Acknowledgments.** The authors are grateful to Reda Messikh for explaining to them the argument described in the Appendix, and to Lincoln Chayes for pointing out that the theory developed in this paper can be used to prove analyticity in $\beta$ of the inverse correlation length.

M. Campanino
Dipartimento di Matematica
Università di Bologna
Piazza di Porta San Donato 5
Bologna
Italy
E-mail: campanin@dm.unibo.it

D. Ioffe
Faculty of Industrial Engineering and Management
Technion
Haifa 3200
Israel
E-mail: ieioffe@ie.technion.ac.il

Y. Velenik
Section de mathématiques
Université de Genève
2-4 rue du Lièvre, CP 64
Geneva
Switzerland
E-mail: Yvan.Velenik@math.unige.ch